\documentclass[11pt]{amsart}
\usepackage{graphicx}
\usepackage{amsmath,amsthm}
\usepackage{amssymb}
\usepackage{pinlabel}
\usepackage[left = 3cm, right = 3 cm , top = 3 cm , bottom = 3cm ]{geometry}
\usepackage[utf8]{inputenc}
\usepackage{cite}
\usepackage{todonotes}
\usepackage{enumerate}

\newtheorem{thm}{Theorem}[section]
\newtheorem{prop}[thm]{Proposition}
\newtheorem{lem}[thm]{Lemma}
\newtheorem{cor}[thm]{Corollary}

\newtheorem{obs}[thm]{Observation}
\newtheorem{question}[thm]{Question}

\newtheorem{notation}[thm]{Notation}
\newtheorem*{eulerclassone}{Euler Class One Conjecture}
\newtheorem*{fullymarkedsurface}{Fully Marked Surface Theorem}
\newtheorem*{maintheorem}{Main Theorem}
\newcommand{\Text}[1]{\text{\textnormal{#1}}}

\theoremstyle{remark}
\newtheorem{remark}[thm]{Remark}

\newtheorem{ex}[thm]{Example}

\theoremstyle{definition}
\newtheorem{definition}[thm]{Definition}

\begin{document}

\title{On Thurston's Euler class one conjecture}
\author{Mehdi Yazdi}
\date{}
\thanks{}
\maketitle

\begin{abstract}
In 1976, Thurston proved that taut foliations on closed hyperbolic 3-manifolds have Euler class of norm at most one, and conjectured that conversely, any integral second cohomology class with norm equal to one is the Euler class of a taut foliation. This is the first from a series of two papers that together give a negative answer to Thurston's conjecture. Here counterexamples have been constructed conditional on the fully marked surface theorem. In the second paper, joint with David Gabai, a proof of the fully marked surface theorem is given.

\end{abstract}

\section{Introduction}

A two-dimensional transversely oriented foliation of a compact, orientable 3-manifold $M$ is called \emph{taut} if every leaf has a closed transversal, where a \emph{transversal} is a closed loop that intersects the leaves of the foliation transversely. Manifolds admitting taut foliations have properties similar to hyperbolic 3-manifolds. In particular, if $M \neq S^2 \times S^1$ then $M$ is \emph{irreducible} \cite{novikov1965topology,rosenberg1968foliations}, i.e. every embedded smooth sphere in $M$ bounds a solid ball. 

For any compact, orientable 3-manifold $M$, Thurston defined natural seminorms on the second homology groups $H_2(M)$ and $H_2(M , \partial M)$, now called the \emph{Thurston norm}\footnote{All (co)homologies have coefficients in $\mathbb{R}$, unless otherwise stated.}. For a compact, orientable surface $S$ with connected components $S_1, \cdots, S_k$, define the \emph{negative part of the Euler characteristic} as 
\[\chi_-(S)=\sum_{i \hspace{1mm}: \hspace{1mm} \chi(S_i)<0} |\chi(S_i)|.\] 
For an integral point $a \in H_2(M)$ or $H_2(M, \partial M)$, the norm of $a$, $x(a)$, is defined as
\[ x(a):= \min \{ \chi_-(S) \hspace{1mm}|\hspace{1mm} [S]=a, \hspace{2mm} S \hspace{1mm} \Text{is properly embedded and oriented} \}. \]
The norm can be extended to rational points linearly and to real points continuously. Up to scaling, Thurston norm is the same as the Gromov's simplicial norm \cite{gabai1983foliations}. Thurston norm on $H_2(M)$ and $H_2(M, \partial M)$ naturally defines dual norms on the dual vector spaces $H^2(M)$ and $H^2(M, \partial M)$. 

A compact, properly embedded, orientable surface $S \subset M$ with no spheres and disks components is called \emph{incompressible} if it has no compressing disk, where an embedded disk $(D,\partial D) \subset (M,S)$ is \emph{compressing} if $D \cap S = \partial D$ and $\partial D$ is homotopically nontrivial in $S$. From now on assume that $M$ is closed and irreducible; see Section \ref{section:euler-class} for manifolds with boundary. Roussarie \cite{roussarie1974plongements} and Thurston \cite{thurston1972foliations, thurston1986norm} showed that taut foliations and incompressible surfaces have an `efficient intersection property'. More precisely, a connected incompressible surface $S$ can be isotoped so that it is either a leaf of the foliation or is transverse to the foliation except at finitely many points of saddle tangencies. In the latter case, the number of tangencies is exactly equal to $|\chi(S)|$ by the Poincar\'{e}--Hopf formula. One can put this into the algebraic language of the \emph{Thurston's inequality}. Let $e(\mathcal{F}) \in H^2(M)$ be the Euler class of the tangent bundle to the foliation. For each embedded, incompressible surface $S$ we have the following inequality for the pairing between the second cohomology and homology groups:
\begin{eqnarray*}
	- \langle e(\mathcal{F}) , [S]\rangle \thinspace \leq  |\chi(S)|.
\end{eqnarray*}
See Section \ref{section:index-sum} for this deduction, due to Thurston. Repeating the same argument for the surface $S$ with the opposite orientation implies that
\begin{eqnarray}
|\langle e(\mathcal{F}) , [S]\rangle| \thinspace \leq  |\chi(S)|.
\label{inequality}
\end{eqnarray} 
Moreover, the incompressibility condition is not required for Inequality (\ref{inequality}) since one can replace $S$ by an incompressible surface $S'$ with $[S] = [S']$ and $|\chi(S')| \leq |\chi(S)|$.

The Euler class of the foliation, $e(\mathcal{F})$, is an element of $H^2(M)$. Therefore, it makes sense to talk about the dual Thurston norm of the Euler class, which we denote by $x^*(e(\mathcal{F}))$. In fact, Inequality (\ref{inequality}) can be written in the abbreviated form
\begin{eqnarray}
x^*(e(\mathcal{F})) \leq 1.
\label{dual-norm}
\end{eqnarray} 

In other words, the Euler class has dual norm at most one. Moreover, the equality happens if $\mathcal{F}$ has any compact leaf of negative Euler characteristic. See Section \ref{section:index-sum}. This puts extreme bounds on the Euler class of a taut foliation. In particular, if $M$ is hyperbolic the number of cohomology classes that can arise as the Euler class of some taut foliation on $M$ is finite. Thurston conjectured that conversely the following happens \cite[Page 129, Conjecture 3]{thurston1986norm}. We call it \emph{the Euler class one conjecture}. A compact 3-manifold is \emph{atoroidal} if every embedded, incompressible torus is boundary parallel.

\begin{eulerclassone}[Thurston--1976] Let $M$ be a closed, orientable, irreducible and atoroidal 3-manifold and let $a \in H^2(M; \mathbb{Z})$ be any integral class with $x^*(a)=1$. Then there is a taut foliation $\mathcal{F}$ on $M$ whose Euler class is equal to $a$.
	\label{conjecture}
\end{eulerclassone}
A compact, orientable, irreducible 3-manifold is \emph{Haken} if it contains an incompressible surface. A compact orientable 3-manifold is \emph{hyperbolic} if its interior admits a complete Riemannian metric of constant sectional curvature $-1$. By Thurston's hyperbolization theorem,  every closed atoroidal Haken 3-manifold is hyperbolic \cite{thurston1986hyperbolic, thurston1998hyperbolic,  thurston1998hyperbolicIII}. See e.g. Otal \cite{otal2001hyperbolization}, and Kapovich \cite{kapovich2001hyperbolic} for the details. Since the manifolds that we consider here have positive first Betti number, they are Haken. Hence $M$ being atoroidal in the statement of the above conjecture is equivalent to $M$ being hyperbolic. Thurston proved Inequality (\ref{inequality}) by showing that there are index sum formulae for both sides of the inequality, and then comparing these sums. See Section \ref{section:index-sum}. It immediately follows from comparing those sums that both sides of (\ref{inequality}) have the same parity when $M$ is closed, i.e. both sides are even since $\chi(S)$ is even. 
\begin{definition}
	Let $M$ be a closed, orientable 3-manifold. An integral class $a \in H^2(M; \mathbb{Z})$ satisfies the \emph{parity condition} if $a \in 2H^2(M; \mathbb{Z})$. An integral class $a \in H^2(M ; \mathbb{R})$ satisfies the \emph{parity condition} if $a = 2b$ for some integral class $b \in H^2(M; \mathbb{R})$.
\end{definition}
In fact, it is known that the integral Euler class of any transversely oriented plane field on a closed, oriented 3-manifold satisfies the parity condition. See Proposition \ref{parity}. Clearly if an integral class $a \in H^2(M ; \mathbb{Z})$ satisfies the parity condition then its image in $H^2(M ; \mathbb{R})$ under the change of coefficients satisfies the parity condition; the converse is also true assuming that $H^2(M ; \mathbb{Z})$ has no $2$-torsion. The parity condition was not explicitly mentioned by Thurston in the Euler class one conjecture, but it was certainly known to him. So we always assume the parity condition as part of the hypotheses of the conjecture. Gabai gave a partial positive answer to the Euler class one conjecture \cite[Page 24, Remark 7.3]{gabai1997problems}; for a proof see Gabai and Yazdi \cite{secondpaper}. By Thurston, the dual unit ball is a convex polyhedron with integral vertices \cite{thurston1986norm}.

\begin{thm}[Gabai]
Let $M$ be a compact oriented irreducible 3-manifold, possibly with toral boundary, and let $a \in H^2(M, \partial M)$ be a vertex of the dual unit ball. Then there is a taut foliation on $M$ whose Euler class is equal to $a$.
\label{gabai}
\end{thm}

A foliation is $C^{\infty,0}$ if the leaves are smoothly immersed with continuous holonomy. We construct counterexamples to the Euler class one conjecture in the general setting of $C^{\infty,0}$ taut foliations, conditional on the fully marked surface theorem. In the early 2000, new techniques were developed that showed many hyperbolic 3-manifolds do not admit taut foliations. See Roberts, Shareshian and Stein \cite{roberts2003infinitely}, and Calegari and Dunfield \cite{calegari2003laminations}. These methods do not apply to our case since Thurston's conjecture is about manifolds with positive first Betti number and so Gabai's theorem \cite{gabai1983foliations} guarantees the existence of taut foliations. Hence, one needs a completely different argument for ruling out taut foliations with a certain Euler class. 

\begin{maintheorem}
There are infinitely many closed hyperbolic 3-manifolds for which the Euler class one conjecture does not hold, i.e. there is some integral second cohomology class with dual Thurston norm equal to one and satisfying the parity condition that is not realized as Euler class of any taut foliation.
\end{maintheorem}

\begin{figure}
\labellist
\pinlabel $(\frac{1}{2},0)$ at 300 55
\pinlabel $(0,\frac{1}{2g-2})$ at 150 135

\pinlabel $(2,0)$ at 470 55
\pinlabel $(0,2g-2)$ at 405 165
\endlabellist

\centering
\includegraphics[width=4 in]{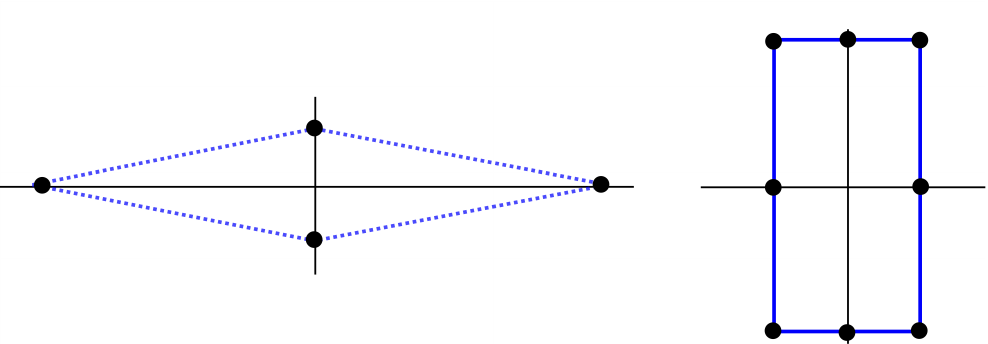}
\caption{Unit ball of the Thurston norm (left), and the dual Thurston norm (right).}
\label{unitballs}
\end{figure}

Our counterexamples are obtained by a suitable Dehn surgery on certain fibered hyperbolic 3-manifolds. \emph{Dehn surgery} is the operation of removing a solid torus from a 3-manifold and gluing it back differently. The constructed counterexamples are explicit in the sense that the monodromy of the fibration map is given in terms of Dehn twists. Moreover, the attaching map corresponding to the Dehn surgery is described. These manifolds are optimal from homological point of view, i.e. their first Betti numbers are equal to two and the unit balls of their Thurston norms have a simple shape (Figure \ref{unitballs}). This is the best that one can hope for, since Gabai's theorem \cite[Theorem 5.5]{gabai1983foliations} implies the truth of the conjecture for 3-manifolds whose first Betti number is equal to one. 

A properly embedded, oriented, possibly disconnected, incompressible surface $S$ is \emph{algebraically fully marked} when in (\ref{inequality}) the equality holds. Note that an algebraically fully marked surface is norm-minimizing, since in this case 
\[ x([S]) \leq |\chi_-(S)| = |\chi(S)| = |\langle  e(\mathcal{F}), [S] \rangle| \leq x([S]), \]
where the last inequality is by (\ref{dual-norm}). The crucial but elementary observation is that any compact leaf is fully marked \cite{thurston1986norm}. See Section \ref{section:index-sum}. The converse, however, cannot be true since one can homotope $\mathcal{F}$ to a new taut foliation without changing the Euler class but with a drastic change in the leaves so that there is no compact leaf any more. The content of the second paper \cite{secondpaper}, joint work with David Gabai, will be a converse to this statement up to homotopy of the plane fields of foliations and under some assumptions. We call this \emph{the fully marked surface theorem}. 

\begin{fullymarkedsurface}[Gabai--Yazdi]
	Let $M$ be a closed hyperbolic 3-manifold, $\mathcal{F}$ be a taut foliation on $M$, and $S$ be an algebraically fully marked surface in $M$. Assume that $S$ is the unique norm-minimizing surface in its homology class up to isotopy. There exists a taut foliation $\mathcal{G}$ on $M$ that has $S$ as a leaf, and such that the oriented plane fields tangent to $\mathcal{F}$ and $\mathcal{G}$ are homotopic.
	\label{fullymarked}
\end{fullymarkedsurface}
We refer the reader to \cite{secondpaper} for a more general statement of the fully marked surface theorem.

\subsection{The key new idea}\ \\

Thurston's Euler class one conjecture predicts that every integral second cohomology class of \underline{dual Thurston norm exactly equal to one} on \underline{closed, hyperbolic 3-manifolds} is realized as the \underline{Euler class} of a taut foliation. The key idea is to examine the conjecture more generally, and ask whether every integral cohomology class of \underline{dual sutured Thurston norm at most one} on \underline{taut sutured 3-manifolds with tori boundary} (not necessarily closed or hyperbolic) can be realized as the \underline{relative Euler class} of a taut foliation. Sutured manifolds were introduced by Gabai for studying taut foliations on 3-manifolds \cite{gabai1983foliations}. One should think about sutures as extra data on the boundary of the 3-manifold indicating how the foliation intersects the boundary. See Subsection \ref{section:sutured-manifold}. In the case of sutured manifolds with toral boundary, by a generalization of the Roussarie--Thurston general position, Inequality (\ref{inequality}) takes the form  
\begin{equation}
|\langle e(\mathcal{F}) , [S]\rangle| \thinspace \leq | \chi_s(S)|,
\label{suturedinequality}
\end{equation} 
where 
\begin{enumerate}
	\item $\chi_s(S)$ is the \emph{sutured Euler characteristic}. See Section \ref{section:sutured-euler-char}.
	\item $e(\mathcal{F})$ is the \emph{relative Euler class}. See Section \ref{section:euler-class} and Definition \ref{def:relative-euler-class}.
\end{enumerate}  

\begin{figure}
	\centering
	\includegraphics[width=3in]{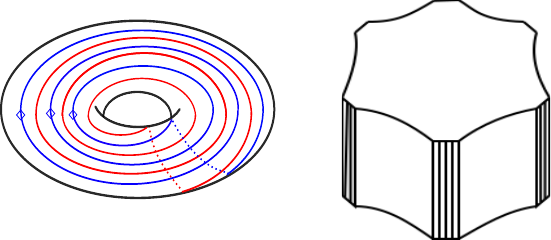}
	\caption{Left: Sutured solid torus $N$. One of the sutures is marked with diamonds in order to make it visually salient. Right: its cross section. }
	\label{solidtorus}
\end{figure}

With this new formulation, one of the first examples to look at is a solid torus $N$ with two parallel sutures on its boundary. Fix a longitude, and assume that each suture goes once around the meridian and three times around the longitude (Figure \ref{solidtorus}). In this case $H_2(N, \partial N)$ is isomorphic to $\mathbb{R}$ and is generated by a meridional disk $D$ of $N$. The sutured Euler characteristic of $D$ is equal to $-2$, and the unit ball and dual unit ball for \emph{sutured Thurston norm} are respectively the intervals $[-\frac{1}{ 2}, \frac{1}{2}]$ and $[-2,2]$ within the real line. See Section \ref{section:sutured-euler-char} for the definition of the sutured Thurston norm.

There are three integral classes in the dual unit ball that satisfy the \emph{relative parity condition}, namely the points $-2$, $0$, and $2$. Here an integral point $a \in H^2(N , \partial N)$ satisfies the relative parity condition if the numbers $\chi_s(D) = -2$ and $\langle a , [D] \rangle$ have the same parity. The two points $\pm 2$ are realized as relative Euler classes of taut foliations comprised of stacks of monkey saddles. See Example \ref{stackexample}. We show in Corollary \ref{nofoliation} that any taut foliation of $N$ is \emph{standard} in the sense that it is obtained from the stack of monkey saddles by a simple operation that does not change the relative Euler class. See Definition \ref{standard}. As a result, the point $0$ is not realized as the relative Euler class of any taut foliation on $N$. This in turn shows that the new formulation of the Euler class one conjecture is violated for the manifold $N$ and the relative Euler class $0$. This gave us some hope that the original Euler class one conjecture also had a chance of not being true.

\subsection{New avenues of research}\ \\

\textbf{Homotopy classification of taut foliations}: Kronheimer and Mrowka proved that on a closed orientable 3-manifold, there are only finitely many homotopy classes of plane fields of taut foliations \cite{kronheimer1997monopoles}. See Gabai \cite{gabai1999essential} for an alternative proof. The tangent bundle of every closed, orientable 3-manifold $M$ is trivial \cite{stiefel1935richtungsfelder, milnor1975characteristic}. Fixing a trivialization $\tau$ of the tangent bundle of $M$, the homotopy class of a transversely oriented plane field $\xi$ is determined by two pieces of information: a cohomology class $\Gamma_\tau \in H^2(M, \mathbb{Z})$ with the property that $2 \Gamma_\tau$ is equal to the Euler class $e(\xi) \in H^2(M ; \mathbb{Z})$, and an element of the affine space $\mathbb{Z}/d(e)\mathbb{Z}$ where $d(e)$ is the \emph{divisibility} of the Euler class $e(\xi)$ modulo torsion. See Gompf \cite[Section 4]{gompf1998handlebody}, Kuperberg \cite[Proposition 2.1]{kuperberg1996}, or Pontrjagin \cite{pontrjagin}. 

The Euler class one conjecture, asking the realizability of extremal cohomology classes, can be understood as a proposed first step in the classification of plane fields of taut foliations on manifolds with positive first Betti number. Although we gave a negative answer to Thurston's conjecture, it is an important question to know exactly which cohomology classes are realized as Euler classes of taut foliations. More generally, given a closed orientable 3-manifold with positive first Betti number, classify homotopy classes of plane fields of taut foliations on the manifold.\\

\textbf{Euler classes of tight contact structures}: Following the work of Bennequin \cite{bennequin1982entrelacements} (or see the English translation \cite{bennequin1989linkings}) and Thurston, Eliashberg proved the analogue of Inequality (\ref{inequality}) for Euler classes of \emph{tight contact structures} \cite{eliashberg1992contact}. Therefore, one would like to know if the Euler class one conjecture has a positive answer for tight contact structures. Every \emph{topological} taut foliation of a closed orientable irreducible 3-manifold can be $C^0$-approximated by a tight contact structure so Euler classes of tight contact structures include Euler classes of taut foliations. See Eliashberg and Thurston \cite{eliashberg1998confoliations}, Kazez and Roberts \cite{kazez2015approximating}, and Bowden \cite{bowden2016approximating}. In particular, there is a chance that the Euler class one conjecture has a positive answer for tight contact structures. 

Colin, Giroux and Honda proved that for every closed oriented 3-manifold there are only finitely many homotopy classes of plane fields that carry tight contact structures \cite{colin2003coarse, colin2009finitude}. Again understanding the set of possible values for the Euler classes of tight contact structures is a proposed first step for the homotopy classification of tight contact structures on closed orientable 3-manifolds with positive first Betti number.\\

\textbf{Euler classes of representations into $\text{Homeo}^+(S^1)$}: Denote by $\text{Homeo}^+(S^1)$ the group of orientation-preserving homeomorphisms of the circle. In what follows, by a representation we mean a representation into $\text{Homeo}^+(S^1)$ unless otherwise stated. Given a taut foliation $\mathcal{F}$ of a closed orientable irreducible atoroidal 3-manifold $M$, the \emph{universal circle construction} of Thurston \cite{thurston1998three}, and Calegari and Dunfield \cite{calegari2006universal} defines a faithful representation 
\[\rho_{\mathcal{F}} \colon   \pi _1(M) \longrightarrow \text{Homeo}^+(S^1). \] 
There is a canonical way to assign an Euler class $e(\rho) \in H^2(M; \mathbb{Z})$ to each such representation using Borel's construction \cite{calegari2006universal}. If $\rho=\rho_{\mathcal{F}}$ is obtained from a taut foliation $\mathcal{F}$, then $e(\rho) = e(\mathcal{F})$; for a proof see Boyer and Hu \cite[Section 6]{boyer2019taut}. Therefore, if the Euler class one conjecture about taut foliations had been true, then for every integral second cohomology class $a$ of dual norm one and satisfying the parity condition there would have been a representation with Euler class equal to $a$. On the other hand, Inequality (\ref{inequality}) holds for Euler classes of representation into $\text{Homeo}^+(S^1)$ according to the \emph{Milnor--Wood inequality} \cite{milnor1958existence, wood1971bundles}. In the case of representations, the Euler class is an element of $H^2(\pi_1(M); \mathbb{Z})$ and is defined as the obstruction to lifting the representation into the \emph{universal central extension} $\widetilde{\text{Homeo}^+(S^1)}$ of $\text{Homeo}^+(S^1)$. Moreover, since $M$ is aspherical, its group cohomology and singular cohomology are naturally isomorphic. Therefore, the set of Euler classes of faithful representations contains the set of Euler classes of taut foliations and it makes sense to ask whether the Euler class one conjecture holds for representations. 

If we ask the same question in one dimension lower that is when $M=S$ is a closed orientable hyperbolic surface, then Milnor \cite{milnor1958existence} and Wood \cite{wood1971bundles} showed that a class $a \in H^2(S; \mathbb{Z})$ is realized as the Euler class of a representation $\rho \colon \pi_1(S) \rightarrow \text{Homeo}^+(S^1)$ if and only if 
\[ |\langle a , [S] \rangle | \leq |\chi(S)|.   \]
The existence of faithful representations of surface groups follows from the work of DeBlois and Kent \cite{deblois2006surface} who proved that all such cohomology classes are realized as Euler classes of faithful representations into the subgroup $\text{PSL}(2, \mathbb{R})$ of $\text{Homeo}^+(S^1)$, answering a conjecture of Goldman \cite{goldman, goldman1988topological}. 

An interesting perspective is to obtain intuition from representations to say something about taut foliations, or vice versa. For example, Boyer, Rolfsen and Wiest have proved that the fundamental group of any closed orientable irreducible 3-manifold with positive first Betti number is \emph{locally indicable} and hence \emph{left-orderable} \cite{boyer2005orderable}, implying that it has a faithful representation into $\text{Homeo}^+(S^1)$ with Euler class zero. Therefore, one can ask whether every closed hyperbolic 3-manifold with $b_1 >0$ admits a taut foliation with Euler class zero.\\

\textbf{Euler classes of pseudo-Anosov flows and quasigeodesic flows}: Inequality (\ref{inequality}) holds for the Euler class of the orthogonal plane bundle to a \emph{pseudo-Anosov flow} (respectively \emph{quasigeodesic flow}) on a closed hyperbolic 3-manifold. See Mosher \cite{mosher1992dynamical,mosher1992dynamicalII}, Calegari and Dunfield \cite{calegari2003laminations}, and Calegari \cite{calegari2006universal}. Therefore, one would like to know whether the Euler class one conjecture holds in this context. 
\subsection{Plan of the proof}\ \\

Section \ref{section:euler-class} gives a careful exposition of the relative Euler class, as well as the index sum. It is shown in Section \ref{section:euler-class} that the index sum is always defined, and it coincides with the evaluation of the relative Euler class for foliations of 3-manifolds with toral boundary. In Section \ref{construction-of-manifolds}, the proposed counterexamples are constructed. We start with a certain fibered hyperbolic 3-manifold $M_f$ with fiber $S$ of genus $g$ and monodromy $f \colon S \rightarrow S$. The final constructed manifold, $M$, is the result of a specific Dehn surgery on $M_f$. The surgery curve is disjoint from a copy of the fiber; in particular, there is a copy of $S$ in $M$ as well. Note that if we cut $M$ along $S$ to get the manifold $M_1 = M \setminus \setminus S$, then $M_1$ is almost a product. 

In Section \ref{section:properties-of-manifolds}, the following properties of the constructed manifold $M$ are established: 
\begin{enumerate}[a)]
	\item The first Betti number of $M$ is equal to $2$, and the unit balls for the Thurston norm and dual Thurston norm of $M$ are as in Figure \ref{unitballs} (Lemma \ref{Thurston-norm}).
	\item $M$ is irreducible and atoroidal (Lemma \ref{atoroidal}). By Thurston's hyperbolization theorem for Haken manifolds, this implies that $M$ is hyperbolic. 
	\item The surface $S \subset M$ is the unique norm-minimizing surface in its homology class up to isotopy (Lemma \ref{unique}). This is proved using the following less known theorem of Gabai regarding norm-minimizing surfaces in a fixed homology class.
	\newtheorem*{disjointsurfaces}{Theorem \ref{disjointsurfaces}}
	\begin{disjointsurfaces}[Gabai \cite{gabai1983foliations}]
		Let $M$ be a closed, orientable 3-manifold. Assume that $P$ and $Q$ are two possibly disconnected, norm-minimizing surfaces in $M$ that are homologous. There is a sequence of possibly disconnected, norm-minimizing surfaces $P=P_0 , P_1 , \cdots , P_n=Q$ with each term in the same homology class as $[P]=[Q]$ such that any two adjacent terms in the sequence can be isotoped to be disjoint in $M$.
	\end{disjointsurfaces}
	
\end{enumerate}
Conditions (2) and (3) above are the hypotheses of the fully marked surface theorem, which will be used in the proof of the main theorem. 

In Section \ref{section:sutured-solid-tori}, we prove that any taut foliation of a sutured solid torus with two sutures is \emph{standard}. See Definition \ref{standard} and Corollary \ref{nofoliation}. In particular, the relative Euler class of the foliation is \emph{maximal} and hence nonzero. See Proposition \ref{maximal euler class} for a statement regarding general sutured solid tori. The proof uses the notion of \emph{based transversal semigroup of a leaf}, which was introduced by Novikov. See Definition \ref{transversal-semigroup}.

We would like to pick a suitable element $a \in 2H^2(M ; \mathbb{Z})$ of dual norm equal to $1$, and show that $a$ is \underline{not} realized as the Euler class of any taut foliation on $M$. We choose the cohomology class $a$ to be the point $(0,2-2g)$ in Figure \ref{unitballs}, right. Note that $\langle a , [S] \rangle = 2-2g = \chi(S)$. Therefore if $a = e(\mathcal{F})$ for some taut foliation $\mathcal{F}$ of $M$, then $S$ is algebraically fully marked with respect to $\mathcal{F}$. By the fully marked surface theorem, if such an $\mathcal{F}$ exists, then we may assume that $S$ is a leaf of $\mathcal{F}$. Therefore we may cut the manifold $M$ and the foliation $\mathcal{F}$ along $S$ and analyze the resulting foliation $\mathcal{F}_1$ on the simpler manifold $M_1$ in order to obtain a contradiction.

By construction, $M_1$ consists of two parts: a product part and a \emph{twisted part}. See Figure \ref{N-and-N'}, right. The twisted part is a solid torus $N'$ with two sutures on its boundary where each suture goes three times in the longitude direction and once in the meridian direction. See Figure \ref{solidtorus}. We show that the initial assumption $e(\mathcal{F})=a$ implies that the restriction of $\mathcal{F}_1$ to $N' \subset M_1$ is a taut foliation with relative Euler class zero. This gives the desired contradiction since any taut foliation of $N'$ is standard and has nonzero relative Euler class (Corollary \ref{nofoliation}).


\subsection{Acknowledgement}

This paper was part of my doctoral thesis at Princeton University. I would like to thank my adviser, David Gabai, for his continuous support and encouragement. Sincere thanks to Yair Minsky who pointed out a mistake in the earlier version and helpful discussions. Many thanks to Shabnam Hossein for helping me with producing the figures. I wish to thank Sergio Fenley, Pierre Dehornoy, and Andr\'{a}s Juh\'{a}sz for helpful conversations and/or comments. Many thanks to the anonymous referee whose constructive comments and corrections improved the paper significantly. The author acknowledges the support by a Glasstone Research Fellowship, and partial support by NSF Grants DMS-1006553 and DMS-1607374.
\section{Preliminaries}

\begin{notation}
	For a metric space $X$ and $A \subset X$, denote the metric completion of $X - A$ with the induced path metric by $X \setminus \setminus A$. Here we are mainly interested in the case that $X$ is a manifold of dimension $2$ or $3$ and $A$ is a submanifold of codimension one. 
	
	For a subset $A$ of $X$, the interior of $A$ is denoted by $A^ \circ$ or $\text{int}(A)$. The number of connected components of $X$ is denoted by $|X|$.
\label{notation}
\end{notation}

\subsection{Foliations}
A two-dimensional \emph{foliation} $\mathcal{F}$ of a 3-manifold $M$ is a partition of $M$ into injectively immersed surfaces that locally looks like the product foliation $\mathbb{R}^2 \times \mathbb{R}$. A \emph{leaf} of the foliation is a connected component of the surfaces in the foliation. The foliation $\mathcal{F}$ is called \emph{transversely oriented} if there is a compatible choice of transverse orientations on the leaves. If the manifold $M$ is oriented, a transverse orientation induces an orientation on each leaf as well. Lickorish \cite{lickorish1965foliation}, Novikov and Zieschang \cite{novikov1965topology} showed that every closed orientable 3-manifold has a foliation, so the existence of a foliation does not give much information about the ambient 3-manifold. However, taut foliations reflect many topological and geometric properties of the ambient manifold. 

For technical reasons, we need to specify the degree of smoothness that we consider here. A foliation $\mathcal{F}$ is called $C^0$, or \emph{topological}, if the holonomy of its leaves is continuous. By Calegari \cite{calegari2001leafwise}, every topological foliation of a 3-manifold is topologically isotopic to a $C^{\infty,0}$ foliation.

\subsection{Suspension foliations}

The exposition for this section is from \cite[Chapter V]{camacho2013geometric}. Let $p \colon E \longrightarrow B$ be a fiber bundle with base $B$,  fiber $F$, and total space $E$. A foliation $\mathcal{F}$ of $E$ is \emph{transverse to the fibers} if 
\begin{enumerate}
	\item Each leaf $L$ of $\mathcal{F}$ is transverse to the fibers and $\dim(L) + \dim(F)= \dim E $.
	\item For each leaf $L$ of $\mathcal{F}$, the restriction map $p \colon L \longrightarrow B$ is a covering map.
\end{enumerate}
Ehresmann showed that when the fiber $F$ is compact, Condition (2) is implied by (1). 

Let $\text{Homeo}(F)$ be the group of homeomorphisms of $F$. Given a fiber bundle and a foliation transverse to the fibers, there is a representation 
\[ \phi \colon \pi_1(B, b_0) \longrightarrow \text{Homeo}(F), \hspace{3mm} b_0 \in B, \]
encoding the \emph{holonomy} of based loops in $B$. Conversely:

\begin{thm}
	Let $B$ and $F$ be connected manifolds. Given a representation 
	\[ \phi \colon \pi_1(B, b_0) \longrightarrow \text{Homeo}(F), \hspace{3mm} b_0 \in B, \]
	there is a fiber bundle $E(\phi)$ over the base $B$ and with fiber $F$, and a foliation $\mathcal{F}(\phi)$ transverse to the fibers of $E(\phi)$ such that the holonomy of $\mathcal{F}(\phi)$ is equal to $\phi$.
	\label{holonomy-rep}
\end{thm}

We are mainly interested in the case that $F = I$ or $S^1$ is one-dimensional, and the image of $\phi$ lies in $\text{Homeo}^+(F)$, that is the group of orientation-preserving homeomorphisms of $F$. The construction is as follows: Let $\tilde{B}$ be the universal cover of $B$. Consider the action of $\pi_1(B,b_0)$ on $\tilde{B} \times F$ defined as 
\[ \gamma \in \pi_1(B , b_0), \hspace{3mm} (\tilde{b} , f ) \in \tilde{B} \times F; \hspace{3mm} \gamma \cdot (\tilde{b} , f) := (\gamma \cdot \tilde{b} \hspace{1mm} , \hspace{1mm} \phi(\gamma) \cdot f),\]
where the action on the first factor is by covering transformations. This action preserves the fibers of $\tilde{B} \times F \xrightarrow{p_1} \tilde{B}$, and induces the structure of a fiber bundle on the quotient 
\[ E(\phi) : = (\tilde{B} \times F) / \pi_1(B,b_0). \]
Moreover, the leaves, $\tilde{B} \times f$, of the product foliation $\tilde{B} \times F$ are preserved by the action of $\pi_1(B, b_0)$. So, this foliation will induce on the quotient space $E(\phi)$ a foliation transverse to the fibers, which we denote by $\mathcal{F}(\phi)$. 

\begin{lem}[\cite{MR1162560}]
	If $F$ is any orientable surface with boundary which is not compact planar and $b$ is a boundary component of $F$, then there are foliations of $F \times I$ ($I$ is a closed interval), transverse to the $I$ factor that have a given holonomy on $b$ and trivial holonomy on all other boundary components. In the remaining case that $F$ is compact planar and not a disk, if $b$ and $b'$ are two boundary components with the induced orientations from $F$, then there exists a foliation transverse to $I$ factor that has a given holonomy $\mu$ on $b$ and $\mu ^{-1}$ on $b'$ and trivial holonomy on all other boundary components.
	\label{transversefoliations}
\end{lem}

The above lemma follows from Theorem \ref{holonomy-rep}, and the fact that every element of $\text{Homeo}^+(I)$ can be written as a commutator. See Li \cite[Lemma 3.1]{li2002laminar}.

\subsection{Haefliger's theorem}
\begin{thm}(Haefliger \cite{haefliger1962varietes})
	Let $\mathcal{F}$ be a codimension-one foliation of a compact $n$-manifold $M$. The union of compact leaves of $\mathcal{F}$ is a compact subset of $M$. Moreover assuming that $\mathcal{F}$ is transversely orientable and $K$ is a compact $(n-1)$-dimensional manifold, the union of leaves of diffeomorphism type $K$ is compact as well.
\end{thm}

%

\subsection{Reeb stability theorem} We need a special case of the Reeb Stability theorem \cite{reeb1952certaines}.

\begin{thm} 
Let $M$ be a compact orientable 3-manifold. Let $\mathcal{F}$ be a transversely oriented codimension-one foliation of $M$ such that $\mathcal{F}$ is transverse to $\partial M$. If some leaf of $\mathcal{F}$ is a sphere (respectively a disk), then $M = S^2 \times S^1$ (respectively $D^2 \times S^1$) with the product foliation. 		
\end{thm}

\subsection{Novikov's and Rosenberg's theorems}
A foliation of a 3-manifold is \emph{Reebless} if it does not contain any \emph{Reeb components}: a foliated solid torus having the boundary torus as a leaf, and with all other leaves being planes spiralling towards the boundary torus \cite{reeb1952certaines, novikov1965topology}. 
\begin{thm}[Novikov--Rosenberg]
 Let $M$ be a compact, orientable 3-manifold. Let $\mathcal{F}$ be a transversely oriented codimension-one foliation of $M$ such that $\mathcal{F}$ has no Reeb components, and $\mathcal{F}$ is transverse to $\partial M$. Then
 \begin{enumerate}
 	\item $M$ is either irreducible or $S^2 \times S^1$ with the product foliation.
 	\item Leaves of $\mathcal{F}$ are $\pi_1$-injective in $M$.
 	\item Every transverse closed curve is homotopically nontrivial.
 \end{enumerate} 
\end{thm}
For (1), Novikov \cite{novikov1965topology} proved that $\pi_2(M)$ is trivial, and Rosenberg \cite{rosenberg1968foliations} showed that $M$ is irreducible. Parts (2) and (3) are due to Novikov \cite{novikov1965topology}. 

\subsection{Thurston Norm} Let $M$ be a compact, orientable 3-manifold. Thurston norm is a seminorm on $H_2(M)$ and $H_2(M, \partial M)$ whose unit ball is a convex polyhedron with rational vertices \cite{thurston1986norm}. See the Introduction for the definition of the Thurston norm.

Associated to any norm $x$ on a vector space $V$ there is a dual norm $x^*$ on the dual vector space $V^*$ defined as 
\begin{eqnarray}
x^*(u) = \sup \{ \langle u , v \rangle \hspace{2mm}|\hspace{2mm} x(v)=1 \} \nonumber.
\end{eqnarray}
This defines the dual Thurston norm on $H^2(M)$ and $H^2(M , \partial M)$.

More generally, for any compact subsurface $A$ of $\partial M$, one can define the Thurston norm on $H_2(M , A)$ along the same lines. See Scharlemann \cite{scharlemann1989sutured} for further generalizations. The corresponding dual Thurston norm is defined on $H^2(M, A) \cong \text{Hom}(H_2(M, A), \mathbb{R})$.

\subsection{Corners}

\begin{figure}
	\labellist
	\pinlabel $3$ at 20 25
	\pinlabel $1$ at 47 40
	\pinlabel $2$ at 65 15
	\pinlabel $4$ at 92 60
	\endlabellist
	
	\centering
	\includegraphics[width= 2 in]{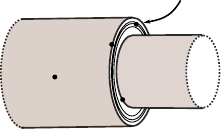}
	\caption{The points 1, 2, 3 and 4 respectively show a convex corner, a concave corner, a tangential point, and a transverse point.}
	\label{corner}
\end{figure}

Consider a two-dimensional foliation of a 3-manifold $M$. Let $p \in \partial M$ be a point. We say that $p$ is a \emph{tangential point} if there is a foliated neighborhood of $p$ in $M$ that is homeomorphic to a foliated neighborhood of $(0,0,0)$ in 
\[ \{ (x,y,z) \hspace{1mm} | \hspace{1mm} x, y \in \mathbb{R}, \hspace{2mm}z \geq 0 \},  \]
where the leaves consist of the planes $z = $ constant. See Figure \ref{corner}. A point $p$ is called a \emph{transverse point} if there is a foliated neighborhood of $p$ in $M$ that is homeomorphic to a foliated neighborhood of $(0,0,0)$ in the foliation of 
\[ \{(x,y,z) \hspace{1mm} | \hspace{1mm} y,z \in \mathbb{R}, \hspace{2mm} x \geq 0\}, \]
where the leaves are the half-planes $z=$ constant.

\begin{notation}
	The transverse (respectively tangential) boundary of $M$ is the closure of the set of transverse (respectively tangential) points in $\partial M$, and is denoted by $\partial_\pitchfork M$ (respectively $\partial_\tau M$).
	\label{notation:boundary}
\end{notation}

A point $p$ is a \emph{convex corner} if there is a foliated neighborhood of $p$ that is homeomorphic to a foliated neighborhood of $(0,0,0)$ in the foliation of 
\[ \{(x,y,z) \hspace{1mm} | \hspace{1mm} x \geq 0 , \hspace{2mm}y \in \mathbb{R}, \hspace{2mm} z \geq 0\}, \]
where the leaves consist of the half-planes $z = $ constant. 

A point $p$ is a \emph{concave corner} if there is a foliated neighborhood of $p$ that is homeomorphic to a foliated neighborhood of  $(0,0,0)$ in the foliation of 
\[ \{  (x,y,z) \hspace{1mm} | \hspace{1mm} y \in \mathbb{R}, \text{ and} \hspace{2mm} \big( x \geq 0  \text{ or } z \geq 0 \big) \}, \]
where the leaves consist of the planes and half-planes $z = $ constant.

\subsection{Sutured manifolds}
\label{section:sutured-manifold}
If $M$ is an oriented 3-manifold, and $S \subset M$ is an embedded oriented surface, then they determine a well-defined transverse orientation for $S$. Conversely a transverse orientation determines an orientation for $S$. 
\begin{definition}
A \emph{sutured manifold} $(M , \gamma)$ is a compact oriented 3-manifold $M$ together with a set $\gamma \subset \partial M$ of pairwise disjoint annuli $A(\gamma)$ and tori $T(\gamma)$. Furthermore, the interior of each component of $A(\gamma)$ contains a \emph{suture}, i.e. a homologically nontrivial oriented simple closed curve. We denote the set of sutures by $s(\gamma)$.

Finally, every component of $R(\gamma) = \partial M - \gamma^\circ $ is oriented. Define $R_+(\gamma)$ (or $R_-(\gamma)$) to be those components of $\partial M - \gamma^\circ$ whose normal vectors point out of (into) $M$. The subsurfaces $R_+(\gamma)$ and $R_-(\gamma)$ are respectively called the \emph{positive} and \emph{negative tangential boundary} of $M$. The orientations on $R(\gamma)$ must be coherent with respect to $s(\gamma)$, i.e. if $\delta$ is a component of $\partial R(\gamma)$ and is given the boundary orientation, then $\delta$ must represent the same homology class in $H_1(\gamma)$ as some suture. 
\end{definition}

\begin{definition}
A transversely oriented foliation $\mathcal{F}$ on a sutured manifold $(M, \gamma)$ is \emph{taut} if
\begin{enumerate}
  	\item $\mathcal{F}$ is tangential to $R(\gamma)$ with the transverse orientation pointing outside (respectively inside) $M$ along $R_+(\gamma)$ (respectively $R_-(\gamma)$). $\mathcal{F}$ is transverse to $\gamma$, and the induced foliation on each component of $\gamma$ is a suspension foliation. 
  	\item Each point on $\partial A(\gamma)$ is a convex corner.
    \item Every leaf of $\mathcal{F}$ has either a closed transversal or a transverse arc starting from $R_-(\gamma) $ and ending on $R_+(\gamma)$. 
\end{enumerate}
$\mathcal{F}$ is \emph{compatible with the sutured structure} if it satisfies Conditions (1) and (2) above. 
\end{definition}

\begin{remark}
	Condition (2) about convex corners is not always mentioned in the definition of a taut foliation of a sutured manifold. However, this is how I learned the concept from Gabai and certainly the foliations constructed by Gabai using sutured manifold hierarchies in \cite{gabai1983foliations} have this property. 
\end{remark}

\subsection{Sutured Euler characteristic and sutured Thurston norm}
\label{section:sutured-euler-char}
Roughly speaking, the sutured Thurston norm in a sutured manifold is defined by doubling the manifold along the tangential boundary (or the transverse boundary), computing the Thurston norm, and then dividing by two. The details are as follows \cite{scharlemann1989sutured, cantwell2013sutured}.

\begin{definition}
Let $(M, \gamma)$ be a sutured 3-manifold. A compact, properly embedded, and oriented surface $S$ in $M$ is \emph{admissible} if for every component $b$ of $\partial S$ either 
\begin{enumerate}
	\item $b$ is an essential simple closed curve in $\gamma$, or
	\item $b \cap A(\gamma)$ is a (possibly empty) union of essential properly embedded arcs. 
\end{enumerate}

\end{definition}
 
\begin{definition}
	Let $S$ be an admissible surface in a sutured manifold $(M, \gamma)$, and $k$ be the total number of arcs in $\partial S \cap A(\gamma)$. Define the \emph{sutured Euler characteristic} of $S$ as 
	\begin{eqnarray}
		\chi_s (S) = \chi (S) - \frac{1}{2} k.
	\end{eqnarray}
\end{definition}

Note that $\chi_s(S)$ depends on the embedding of $S$ and can change under an isotopy.  

\begin{remark}
Let $\mathrm{T}S$ be the tangent bundle of $S$. For an admissible  surface $S$, there is a non-vanishing section of $\mathrm{T} S |\partial S$ canonical up to homotopy: pick a Riemannian metric on $S$ and consider a section of unit length that is tangential to $\partial S$ along $\partial S \cap R(\gamma)$, and that is pointing inward/outward along components of $\partial S \cap \gamma$ in an alternate fashion. See Figure \ref{normal-vectors}, right.

By Poincar\'{e}--Hopf formula, the Euler characteristic of a compact surface $S$ is the obstruction to the existence of a non-vanishing section of the tangent bundle of the surface that always points outward (or always inward, or always tangential) along each boundary component. Given a vector field on $S$ with (generalized) Morse singularities that is entirely tangential or entirely transverse to each component of $\partial S$, the Euler characteristic can be computed as sum of the indices of the singularities. Here the index of a saddle point (respectively \emph{generalized saddle point}) is $-1$ (respectively $1-\frac{n}{2}$ where $n$ is the number of \emph{prongs}) and the index of a center point is $+1$. This picture generalizes to an admissible surface $S$ in a sutured manifold; in this case the sutured Euler characteristic is the obstruction to the existence of a non-vanishing section of $\mathrm{T} S$ whose restriction to $\partial S$ coincides with the section coming from the sutured structure of the manifold as described above. 
\label{sutured-euler-char-poincare-hopf}
\end{remark}

The following is a natural generalization of the Thurston norm to sutured manifolds \cite{scharlemann1989sutured, cantwell2013sutured}. 
\begin{definition}
	Let $(M,\gamma)$ be a sutured manifold, and $S$ be an admissible surface. If $S_1 , \cdots, S_k$ are the components of $S$, define 
	\begin{eqnarray}
	 \chi_-^s(S) = \sum_{i \hspace{1mm}: \hspace{1mm} \chi_s(S_i)<0} |\chi_s(S_i)|.
	\end{eqnarray} 
	Define the \emph{sutured Thurston norm} of an integral class $a \in H_2(M,\partial M)$ by minimizing $\chi_-^s(S)$ over all admissible representatives $S$ of $a$. Define the sutured Thurston norm on $H_2(M , \partial M)$ by extending to rational points linearly, and to real points continuously.
\end{definition}

\begin{remark}
In Cantwell and Conlon \cite{cantwell2013sutured}, the quantity $\chi_-^s(S)$ is defined for more general surfaces; in particular inessential arcs are allowed in $S \cap A(\gamma)$. However, considering the class of admissible surfaces results in the same definition for the sutured Thurston norm.
\end{remark}

\begin{ex}
	Let $M$ be a sutured solid torus with two sutures. Pick a longitude and assume that each suture goes twice around the longitude and once around the meridian. Let $F$ be the meridional disk of $M$. Then $F$ intersects the sutures in $k=4$ arcs and $\chi_s(F) = 1 - \frac{1}{2} \times 4 = -1$. The flow lines of a sample vector field on $F$ are drawn in Figure \ref{stackofchairs}, center. Note that the unit tangent vectors to the flow lines are compatible with the section on $\partial F$ coming from the sutured  structure of $M$. There is only one Morse singularity, and it is of saddle type; hence the index sum gives the same number $-1$ for the sutured Euler characteristic.  
	
	The homology group $H_2(M, \partial M)$ is generated by the meridional disk $F$. Hence the unit balls for the sutured Thurston norm and dual norm are the intervals $[-1,1]$ and $[-1,1]$.
	\label{Ex-stack-of-chairs}
\end{ex}

\subsection{Pseudo-Anosov maps} 
Let $S := S_g$ be a closed, orientable surface of genus $g$. A \emph{multicurve} is a union of distinct (up to isotopy) and disjoint essential simple closed curves on $S$. Given a homeomorphism $\phi$ of $S$, one can look at the action of $\phi$ on the set of multicurves on $S$. Pseudo-Anosov homeomorphisms can be characterized by the property that they do not fix the isotopy class of any multicurve \cite{thurston1988geometry}. See Fathi, Laudenbach and Po\'{e}naru \cite{fathi1979travaux} or the English translation \cite{fathi2012thurston} for an exposition of the rich theory of pseudo-Anosov maps.

Thurston's hyperbolization theorem for fibered 3-manifolds states that the mapping torus of a homeomorphism $\phi \colon S \longrightarrow S$ is a hyperbolic manifold if and only if $\phi$ is isotopic to a pseudo-Anosov map \cite{thurston1998hyperbolic}. See Otal \cite{otal2001hyperbolization}. 

\subsection{Penner's construction of pseudo-Anosov maps}
Thurston gave the first hands-on construction of pseudo-Anosov maps in terms of Dehn twists \cite{thurston1988geometry}. His construction made use of twists along two curves $\alpha$ and $\beta$ such that $\alpha \cup \beta$ fills the surface. Penner gave another construction using opposite twists along multicurves. Let $\alpha = a_1 \cup \dots \cup a_m$ and $\beta = b_1 \cup \dots \cup b_n$ be two multicurves on $S$ such that $\alpha \cup \beta$ \emph{fills} the surface, meaning that each component of $S - \alpha \cup \beta$ is topologically a disk. Let $\tau_{a_i}$ be the positive (right handed) \emph{Dehn twist} along $a_i$. Define $\tau_{b_j}$ similarly. Penner's theorem states that any word in $\tau_{a_i}$ and $\tau_{b_j}^{-1}$ is pseudo-Anosov provided that each $\tau_{a_i}$ and $\tau_{b_j}^{-1}$ is used at least once \cite{penner1988construction}. Note that we are doing positive Dehn twist along the curves in one multicurve and negative twists along the other one.

\begin{figure}
	\labellist
	\pinlabel $a_1$ at 18 58
	\pinlabel $a_2$ at 115 58
	\pinlabel $a_3$ at 215 58
	\pinlabel $b_1$ at 55 62
	\pinlabel $b_2$ at 180 62 
	\endlabellist
	\centering
	\includegraphics[width= 2 in]{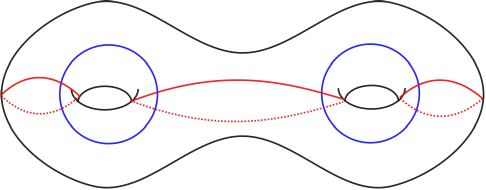}
	\caption{Two multicurves that together fill the surface}
	\label{pennercurvesexample}
\end{figure}

\begin{ex}
	Let $S$ be a closed orientable surface of genus two. Define the multicurves $\alpha = a_1 \cup a_2 \cup a_3$ and $\beta = b_1 \cup b_2$ to be as in Figure \ref{pennercurvesexample}. It can be easily seen that $\alpha \cup \beta$ fills the surface. By Penner's theorem, the map $f = \tau_{a_1}^2 \circ \tau_{a_2} \circ \tau_{b_2}^{-3} \circ \tau_{a_3} \circ \tau_{b_1}^{-1} \circ \tau_{a_1}$ is pseudo-Anosov. 
\end{ex}

\subsection{Gromov simplicial norm} Let $X$ be a compact, oriented manifold possibly with boundary, and let $C_*(X)$ and $C_*(X , \partial X)$ respectively be the real chain complex and the relative real chain complex of $X$. Each $c \in C_*(X)$ (respectively $C_*(X, \partial X)$) is a finite linear combination of singular simplices in $X$ (which do not lie in $\partial X$, and) with real coefficients, i.e. $c = \sum_i \hspace{1mm} r_i \sigma_i$ where $r_i \in \mathbb{R}$ and $\sigma_i$ are singular simplices. Define the $\ell^1$-norm of $c$ as 
\[ ||c|| = \sum_i \hspace{1mm} |r_i|. \]
This induces a seminorm on $H_*(X)$ and $H_*(X, \partial X)$, called the \emph{Gromov norm}, by setting the norm of $\alpha \in H_*(X)$ or $ H_*(X,\partial X)$ to be 
\[ g(\alpha) := \inf ||z||, \]
where $z$ varies between all cycles or relative cycles representing $\alpha$. See Gromov \cite{gromov1982volume}. It follows from the definition that the Gromov norm is decreasing under pushforward.  

In \cite{gabai1983foliations}, Gabai proved the following two theorems, both of which were previously conjectured by Thurston \cite{thurston1986norm}. 
\begin{thm}[Gabai]
	Let $M$ be a compact, oriented 3-manifold. Then on $H_2(M)$ or $H_2(M, \partial M)$, $x = \frac{1}{2}g$ where $x$ denotes the Thurston norm, and $g$ denotes the Gromov norm.
	\label{proportional}
\end{thm}

\begin{thm}[Gabai]
	Let $M$ be a compact, oriented 3-manifold. Let $p \colon \tilde{M} \rightarrow M$ be an $n$-fold covering map and let $z \in H_2(M) = H^1(M, \partial M)$ or $z \in H_2(M, \partial M) = H^1(M)$. Then $n(x(z)) = x(p^*(z))$ where $x(z)$ denotes the Thurston norm of $z$. 
	\label{gabai-covering}
\end{thm}

The next proposition describes the behavior of the dual Thurston norm under (covering) maps between 3-manifolds, making use of the above two theorems of Gabai. 
\begin{prop}
	Let $M$ and $N$ be compact, orientable 3-manifolds. Let $p \colon N \longrightarrow M$ be a map sending $\partial N$ to $\partial M$, and $p^* \colon H^2(M, \partial M) \longrightarrow H^2(N , \partial N)$ be the induced map on real cohomology equipped with the dual Thurston norm. Then 
	\begin{enumerate}
		\item The map $p^*$ is norm decreasing. 
		\item The map $p^*$ preserves the norm if $p$ is a covering map. 
	\end{enumerate}
	\label{covering}
\end{prop}

\begin{proof}
	Let $a \in H^2(M, \partial M)$ be a point with dual norm equal to $k$. For each $h \in H_2(N, \partial N)$ we have:
	\begin{eqnarray}
	\langle p^*(a) , h \rangle = \langle a , p_*(h) \rangle \leq  k \cdot x\big(p_*(h)\big) \leq k \cdot x(h). \nonumber
	\end{eqnarray}
	The last inequality is true since the Thurston norm and the Gromov norm are proportional (Theorem \ref{proportional}) and the Gromov norm is decreasing under pushforward. This implies that $p^*(a)$ has dual norm at most $k$, so $p^*$ is norm decreasing. \\
	
	Now suppose $p$ is a covering map. Let $a \in H^2(M, \partial M)$ be a point with dual norm equal to $k$. There is an integral point $[F] \in H_2(M, \partial M)$ represented by a properly embedded oriented surface $F$ such that 
	\begin{eqnarray}
	\langle a , [F] \rangle = k \cdot x([F]) \nonumber.
	\end{eqnarray}
	This might need some explanation: Let $B_x$ be the unit ball of the Thurston norm on $H_2(M,\partial M)$. The point $a$ is a linear functional on $H_2(M, \partial M)$, and so the supremum of $\{ \langle a , v \rangle \hspace{1mm}| \hspace{1mm} v \in B_x \} $ happens at a vertex of the polyhedron $B_x$, which is a rational point. By scaling this vertex, we obtain an integral point $[F]$.\\
	
	Let $f =[p^{-1}(F)] \in H_2(N , \partial N)$. Then 
	\begin{eqnarray}
	\langle p^*(a) , f \rangle = \langle a , p_*(f) \rangle = \langle a , \mathrm{deg}(p) \cdot [F] \rangle = k \cdot \mathrm{deg}(p) \cdot x\big([F]\big) = k \cdot x(f) \nonumber.
	\end{eqnarray}
	Here the last equality, $x(f) = \mathrm{deg}(p) \cdot x\big([F]\big)$, holds by Theorem \ref{gabai-covering} since $p$ is a covering map. 
\end{proof}

%
%
%

\subsection{Incompressible and $\partial$-incompressible surfaces} Let $M$ be a compact, orientable 3-manifold, and $S$ be a compact, properly embedded, orientable surface with no sphere and disk components. A \emph{$\partial$-compressing disk} for $S$ is an embedded disk $D$ such that $\partial D = \alpha \cup \beta$ with $\alpha \subset S$ an essential arc and $\beta \subset  \partial M$ an arc such that $\alpha \cap \beta = \partial \alpha=\partial \beta$ and $D \cap S = \alpha $. Here $\alpha \subset S$ being \emph{essential} means that there is no embedded disk $D' \subset S$ with $\partial D' = \alpha \cup \beta'$ where $\beta' \subset \partial S$. The surface $S \subset M$ is \emph{$\partial$-incompressible} if it does not admit any $\partial$-compressing disk. We need the following two basic facts about incompressible and $\partial$-incompressible surfaces.

\begin{lem}
	Every connected incompressible surface in a solid torus is a $\partial$-parallel annulus.
	\label{incompressible-surface-in-torus}
\end{lem}
For the proof of the above lemma see Martelli \cite[Proposition 9.3.16]{martelli2016introduction}.

\begin{lem}
	Let $K$ be a connected, compact, orientable surface, and $F \subset K \times [0,1]$ be an incompressible surface with $\partial F \subset \partial K \times (0,1)$. Assume further that $F$ does not admit any $\partial$-compressing disk $D$ with $\partial D = \alpha \cup \beta$, where $D \cap F = \alpha$ is an essential arc, $\beta \subset \partial K \times [0,1]$ is an arc, and $\partial \alpha = \partial \beta = \alpha \cap \beta$. Then each component of $F$ is isotopic to a horizontal surface, i.e. $K \times \{ t \}$ for some $t \in [0,1]$.
	\label{incompressible-product}
\end{lem} 


For the proof of the above lemma for closed surfaces see Martelli \cite[Proposition 9.3.18]{martelli2016introduction}; the proof works for surfaces with boundary without much modification.

\section{Relative Euler class} 
\label{section:euler-class}

Given a plane bundle over base $K$, and a subcomplex $L \subset K$, the definition of a relative Euler class as an element of $H^2(K, L)$ depends on the choice of a non-vanishing section on $L$. In \cite{thurston1986norm}, Thurston considered relative Euler classes of taut foliations on 3-manifolds with boundary; however, he did not explicitly mention a choice of such section. In fact, when the 3-manifold has a surface of genus $>1$ in its tangential boundary, such a non-vanishing section does not exist over the entire boundary for Euler characteristic reasons. Therefore, some care is needed in understanding Corollary 1 of \cite[Page 118]{thurston1986norm} correctly. In Thurston's notation, $\chi$ denotes the Euler class.

\begin{thm}[Thurston - Corollary 1 of \cite{thurston1986norm}]
Let $M$ be an oriented 3-manifold and $\mathcal{F}$ a codimension-one, transversely oriented foliation of $M$. Suppose that $\mathcal{F}$ contains no Reeb components and each component of $\partial M$ is either a leaf of $\mathcal{F}$ or a surface $T$ such that $\mathcal{F}$ is transverse to $T$ and each leaf of $\mathcal{F}$ which intersects $T$ also intersects a closed transverse curve (e.g. $\mathcal{F}|T$ has no Reeb components).  
\[ x^*(\chi(\mathrm{T} \mathcal{F})) \leq 1 \]
holds 
\begin{enumerate}[a)]
	\item in $H^2(M)$,
	\item in $H^2(M, \partial M)$.
\end{enumerate}
\end{thm}

Therefore, we think that it is beneficial to give a detailed account of the relative Euler class of foliations of 3-manifolds; although we do not claim originality for the results of this section. Thurston proved that the evaluation of the relative Euler class on a surface $S$ is equal to an index sum; see Equation (\ref{index-sum}) in Section \ref{section:index-sum} and the proposition in \cite[Page 116]{thurston1986norm}. In this section, we start with the index sum as our definition and show that the index sum determines a well-defined Euler class in $H^2(M, \partial_\pitchfork M)$, but not necessarily in $H^2(M, \partial M)$. See Notation \ref{notation:boundary}, and Remark \ref{higher-genus-boundary}. We also show that when $\partial M$ is a union of tori, the index sum determines a well-defined Euler class in $H^2(M, \partial M)$. Therefore, probably Thurston meant the Euler class as an element of $H^2(M, \partial_\pitchfork M)$ in Corollary 1 of \cite{thurston1986norm}. 

\subsection{Obstruction theory} In this subsection, we define the Euler class from the viewpoint of obstruction theory. We refer the reader to \cite{ steenrod1999topology, milnor1975characteristic
} or \cite[Section 4]{candel2003foliations} for further details. 

Let $K$ be a finite CW complex and $L \subset K$ be a subcomplex. Denote the $n$-skeleton of $K$ by $K^n$. Let $\mathfrak{B}$ be an oriented circle bundle over $K$. Fix a section $\mathfrak{s}$ of $\mathfrak{B}$ over $L$. We are interested in successive extensions of $\mathfrak{s}$ over $K^n \cup L$ as we increase $n = 0 ,1 , 2 , \cdots$. Hence when $L = \emptyset$, this specializes to the problem of finding a section of $\mathfrak{B}$ over $K^n$.

The section $\mathfrak{s}$ can always be extended to $K^0 \cup L$ by picking arbitrary values on the $0$-cells of $K$. Moreover, $\mathfrak{s}$ can be extended to $K^1$ since the fiber $S^1$ is connected. Let $\mathfrak{s}_1$ be an extension of $\mathfrak{s}$ over $K^1 \cup L$; in general $\mathfrak{s}_1$ is not unique even up to homotopy. The first obstruction can happen when we would like to extend $\mathfrak{s}_1$ over $K^2 \cup L$. The obstruction to the existence of an extension of $\mathfrak{s}_1$ to $K^2 \cup L$ can be measured by a cochain $c(\mathfrak{B}, \mathfrak{s}_1)$, called the \emph{obstruction cochain}. If we denote the relative $2$-chains for $(K, L)$ by $C_2(K, L) := C_2(K)/C_2(L)$, then the cochain $c$ assigns integer values to elements of $C_2(K , L)$ in the following way: Given a $2$-cell $\sigma \colon \Delta^2 \rightarrow K$, pick a trivialization of $\mathfrak{B}$ over $\Delta^2$ to identify $\mathfrak{B}|\Delta^2 \cong \Delta^2 \times S^1$. Under this identification, the restriction of $\mathfrak{s}_1$ to $\partial \Delta^2$ determines a map $\partial \Delta^2 \rightarrow S^1$. By definition, the cochain $c$ assigns the degree of the map $\partial \Delta^2 \rightarrow S^1$ to $\Delta^2$. This cochain is in fact a relative cocycle.\\

The obstruction cochain $c$ in general depends not only on $\mathfrak{s}$ but also on the chosen extension $\mathfrak{s}_1$ to the 1-skeleton. However, as we vary $\mathfrak{s}_1$, the cohomology class of $c(\mathfrak{B}, \mathfrak{s}_1)$ remains fixed inside $H^2(K , L ; \mathbb{Z})$; this cohomology class is called the \emph{Euler class of $\mathfrak{B}$ relative to $\mathfrak{s}$} and is denoted by $e(\mathfrak{B}, \mathfrak{s})$. The Euler class $e(\mathfrak{B}, \mathfrak{s})$ remains invariant under a homotopy of $\mathfrak{s}$, and it vanishes if and only if $\mathfrak{s}$ can be extended to a section over $K^2 \cup L$. If $\mathfrak{s}$ can be extended to a section over $K^2 \cup L$, then $\mathfrak{s}$ can be extended over the entire base $K$. To see this, assume that $n \geq 2$ and we would like to extend a given section $\mathfrak{s}_n$ over $K^n \cup L$ to a section over $K^{n+1} \cup L$. Let $\sigma \colon \Delta^{n+1} \rightarrow K$ be an $(n+1)$-cell and choose a trivialization of $\mathfrak{B}$ over $\Delta^{n+1}$. Using the trivialization $\mathfrak{B}|\Delta^{n+1} \cong \Delta^{n+1} \times S^1$, the restriction of $\mathfrak{s}_n$ to $\partial \Delta^{n+1}$ determines a map $\partial \Delta^{n+1} \rightarrow S^1$ which is homotopically trivial since $\pi_n(S^1)$ is trivial for $n \geq 2$. Therefore for each $\Delta^{n+1}$, the map $\partial \Delta^{n+1} \rightarrow S^1$ can be extended to a map $\Delta^{n+1} \rightarrow S^1$, i.e. $\mathfrak{s}_n$ can be extended to a section over $K^{n+1} \cup L$. Hence, the Euler class $e(\mathfrak{B}, \mathfrak{s})$ is the obstruction to extending $\mathfrak{s}$ to a section over the entire base $K$ as well.

Now assume that $\mathfrak{B}$ is an oriented plane bundle over $K$, and equip each fiber with an inner product where the inner product varies continuously. This is possible using a partition of unity and the fact that any convex combination of positive definite inner products is again such an inner product. Then the Euler class of $\mathfrak{B}$ is defined as the Euler class of the associated oriented unit circle bundle. 

\subsection{Euler class of a foliation relative to the transverse boundary}
\label{section:rel-euler-class-transverse-boundary}

Let $M$ be a compact oriented 3-manifold, and $\mathcal{F}$ be a transversely oriented foliation on $M$ such that each component of $\partial M$ is either a leaf of $\mathcal{F}$ or is transverse to $\mathcal{F}$. Define the relative Euler class $e(\mathrm{T} \mathcal{F}, \mathfrak{s}) \in H^2(M , \partial_\pitchfork M; \mathbb{Z})$ using the outward pointing section $\mathfrak{s}$ of $\mathrm{T} \mathcal{F}$ along $\partial_\pitchfork M$. We can alternatively use the inward pointing section of $\mathrm{T}\mathcal{F}|\partial_{\pitchfork} M$, or use the section that is tangential to $\mathcal{F}|\partial_{\pitchfork}M$ (with either orientation) for defining the relative Euler class. This is because all the mentioned sections are homotopic through non-vanishing sections of $\mathrm{T} \mathcal{F}|\partial_\pitchfork M$; the homotopy consists of simultaneous rotation in the planes of $\mathrm{T}\mathcal{F}|\partial_{\pitchfork} M$. We refer to the relative Euler class by $e(\mathcal{F})$ if the choice of $\mathfrak{s}$ is clear.

\subsection{Index sum, and relative Euler class}
\label{section:index-sum}

\begin{figure}
	\centering
	\begin{minipage}{6 cm}
		\labellist
		\pinlabel $N(x)$ at 200 195
		\pinlabel $n(x)$ at 200 110
		\pinlabel $p(\nu_S(x))$ at 100 130
		\pinlabel $\nu_S(x)$ at 133 192 
		\pinlabel $r(x)$ at 145 85
		\pinlabel $T_x(S)$ at 295 195
		\pinlabel $T_x(\mathcal{F})$ at 36 230
		\endlabellist
		
		\centering
		\includegraphics[width=2.3 in]{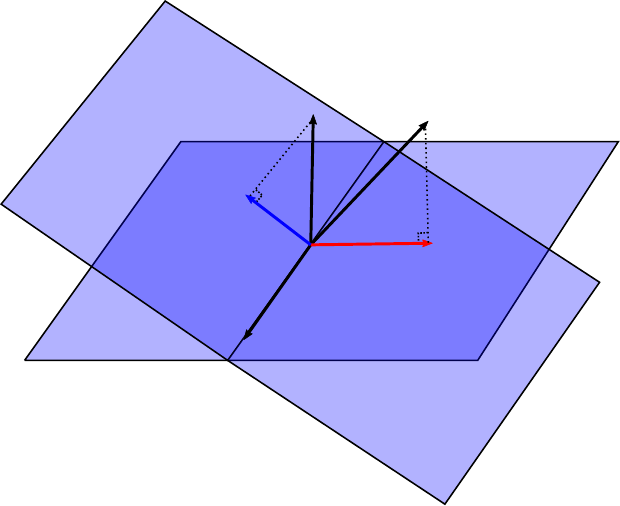}
	\end{minipage}
	\begin{minipage}{6 cm}
		\labellist
		\pinlabel $\nu_S$ at 90 50
		\pinlabel $r(x)$ at 85 72
		\pinlabel $R_-(\gamma)$ at 50 65
		\endlabellist
		
		\centering
		\includegraphics[width =2.4 in]{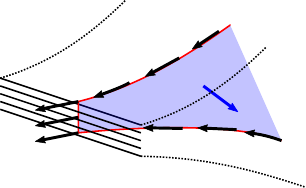}
	\end{minipage}
	
	\caption{Left: the relative position of the tangent plane to the surface $S$ and the foliation $\mathcal{F}$, and various associated vectors. Right: the local picture of an admissible surface intersecting an annulus suture, and the vector field $r$. Note there is another possibility that $\nu_S$ points in the opposite direction, in which case all the arrows for the vector field $r$ are reversed.
	}
	\label{normal-vectors}
\end{figure}

\begin{definition}
	Let $M$ be a sutured manifold, $\mathcal{F}$ be a foliation of $M$ that is compatible with the sutured structure, and $S$ be an admissible surface. Isotope $S$, keeping $\partial S \cap \partial \gamma$ fixed, such that each component of $\partial S \cap \gamma$ is either transverse to $\mathcal{F}|\gamma$ or is a leaf of $\mathcal{F}| \gamma$, and furthermore $\mathcal{F}$ is transverse to $S$ in a neighborhood of $\partial S$. This is possible since $\mathcal{F}|\gamma$ has no Reeb annuli.
	
	Pick a Riemannian metric on $M$. For $x \in S$, let $N(x)$ be the oriented unit normal vector to $\mathcal{F}$, and $n(x)$ be the projection of $N(x)$ onto the tangent plane $\mathrm{T}_x(S)$. Then $n(x)$ is non-vanishing along $\partial S$. Denote the tangent bundle of $\mathcal{F}$ by $\mathrm{T}\mathcal{F}$, and the unit tangent bundle by $\lambda$. Let $r(x) \in \mathrm{T}_x(\mathcal{F})$ be the unit tangent vector to the induced singular foliation $\mathcal{F}|S$ satisfying the condition that $(r(x),n(x))$ is a basis giving the orientation of $\mathrm{T}_x(S)$. See Figure \ref{normal-vectors}, where $\nu_S$ denotes the positive unit normal vector to $S$. Then $r(x)$ is a section of $\lambda$ defined on $S$ minus singular points of $\mathcal{F}|S$. In particular, $r|\partial S$ is a non-vanishing section of $\mathrm{T}\mathcal{F}|\partial S$. See Figure \ref{normal-vectors}, right. 
	Define $\text{Ind}(\mathcal{F}, S)$ as the obstruction to extending the section $r|\partial S$ to a non-vanishing section of $\mathrm{T} \mathcal{F}|S$, i.e. if $[S] \in H_2(S, \partial S)$ is the fundamental homology class of $S$ then 
	\begin{eqnarray*}
		\text{Ind}(\mathcal{F}, S) : = e(\mathrm{T} \mathcal{F}|S , r|\partial S) \cap [S] \in \mathbb{Z}.
	\end{eqnarray*}
	
	\label{def:ind}
\end{definition}

\begin{remark}
$\text{Ind}(\mathcal{F}, S)$ is well-defined, and is invariant under an isotopy of $S$ relative to its boundary. To see this assume that $S$ is admissible, and isotope $S$, keeping $\partial S \cap \partial \gamma$ fixed, to respectively $S_1$ and $S_2$ such that each component of $\partial S_i \cap \gamma$ is either a leaf of $\mathcal{F}|\gamma$ or is transverse to $\mathcal{F}|\gamma$, and $\mathcal{F}$ is transverse to $S_i$ in a neighborhood of $\partial S_i$.

Then $S_1$ is isotopic to $S_2$ through surfaces satisfying the above boundary condition. Let $F \colon S \times [0,1] \rightarrow M$ be such an isotopy with $F|(S \times 0) = S_1$ and $F|(S \times 1) = S_2$. By pulling back the bundle $\mathrm{T} \mathcal{F}$ using the map $F$, we obtain an isomorphism between the bundles $\mathrm{T} \mathcal{F}|S_i$ for $i = 1, 2$ that sends the corresponding sections on $\partial S_i$ to each other. Therefore, the obstruction numbers associated to them are equal as well. 
\end{remark}

Following Thurston \cite{thurston1986norm}, there is an index sum formula for $\text{Ind}(\mathcal{F}, S)$. Isotope $S$, keeping it fixed in a neighborhood of $\partial S$, such that $S$ and $\mathcal{F}$ are transverse to each other except at a finite number of tangencies that are of \emph{saddle} or \emph{center} type. Then
\begin{eqnarray}
\text{Ind}(\mathcal{F}, S) = \sum_{\substack{\Text{tangent }\\ \Text{points }p}}  \text{sign}(p) \cdot  i(p).
\label{index-sum}
\end{eqnarray}
Here $i(p)$ is the \emph{index} of the tangency point. Saddle points have index $-1$ and center points have index $+1$. By definition, $\text{sign}(p) \in \{ -1 , +1\}$ is equal to $+1$ exactly when the transverse orientations of $\mathrm{T}_p(S)$ and $\mathrm{T}_p(\mathcal{F})$ agree. 

We briefly explain why this index sum formula holds; for details see Thurston \cite[Pages 115--118]{thurston1986norm}. The number $\text{Ind}(\mathcal{F}, S)$ is the obstruction to extending the section $r|\partial S$ to a section of $\lambda$ over $S$. By construction, $r(x)$ is a section of $\lambda$ defined on $S$ minus finitely many singular points of $\mathcal{F}|S$. For each such point $p$, we can calculate the degree of the map $r_p \colon \partial D_p \rightarrow S^1$ where $D_p$ is a disk neighborhood of $p$ in $S$ inducing an orientation on $\partial D_p$, and $r_p$ is the restriction of $r$ to $\partial D_p$ with $D_p \times S^1$ the local trivialization of $\lambda$. By summing the degrees over such points $p$, the number $\textrm{Ind}(\mathcal{F}, S)$ is calculated and gives the formula (\ref{index-sum}). The term $\text{sign}(p)$ in the formula comes from the fact that in computing the degree of the map $r_p \colon \partial D_p \rightarrow S^1$, $\partial D_p$ is oriented as the boundary of $D_p$ (or equivalently using the transverse orientation of $S$), whereas $S^1$ is oriented using the transverse orientation $N(x)$ of the foliation \cite[Figure 5]{thurston1986norm}. 

\begin{remark}
Thurston proved the index sum formula in a more general setting where \emph{circle tangencies} between $S$ and $\mathcal{F}$ are allowed as well. A neighborhood of the circle tangency might look like the graph of $1-(x-5)^2 \hspace{4mm}  4 \leq x \leq 6$ revolved about the $z$-axis. Circle tangencies do not contribute to the index sum. See Thurston \cite[Pages 115--118]{thurston1986norm}.
\end{remark}

The Poincar\'{e}--Hopf formula for $\mathcal{F}|S$ shows that the same index sum without the term $\text{sign}(p)$ is equal to the sutured Euler characteristic of $S$. See Remark \ref{sutured-euler-char-poincare-hopf}.
\begin{eqnarray}
\chi_s(S) = \sum_{\substack{\Text{tangent }\\ \Text{points }p}} i(p).
\end{eqnarray}
Thurston observed that if $S$ is a compact leaf of $\mathcal{F}$, then 
\begin{eqnarray}
\text{Ind}(\mathcal{F}, S) = \pm \hspace{1mm} \chi (S) ,
\label{compact-leaf}
\end{eqnarray} 
where the $\pm$ sign is equal to $+1$ exactly when the transverse orientations of $S$ and $\mathcal{F}$ agree \cite{thurston1986norm}. This is because when $S$ is a leaf, the restriction of the tangent bundle $\mathrm{T}\mathcal{F}$ to $S$ coincides with the tangent bundle of $S$, up to orientation.

\begin{prop}
	Let $M$ be a sutured 3-manifold, $\mathcal{F}$ be a foliation of $M$ that is compatible with the sutured structure. Then $\text{Ind}(\mathcal{F}, \cdot) $ determines a cohomology class in $H^2(M , \partial_\pitchfork M ; \mathbb{R}) \cong \text{Hom} (H_2(M , \partial_\pitchfork M), \mathbb{R})$ that coincides with the real Euler class relative to the outward pointing section along $\partial_\pitchfork M$.

\end{prop}

\begin{proof}
Let $\mathfrak{s}_0$ be the unit outward pointing section of $\mathrm{T} \mathcal{F}|\partial_\pitchfork M$, and $e(\mathrm{T}\mathcal{F} , \mathfrak{s}_0) \in H^2(M, \partial_{\pitchfork} M)$ be the Euler class of $\mathrm{T} \mathcal{F}$ relative to $\mathfrak{s}_0$. Let $[S] \in H_2(M, \partial_\pitchfork M)$ be a homology class. There is an admissible surface $S$ whose homology class is $[S]$. Using the naturality property of the Euler class for the map corresponding to the embedding of $S$ in $M$, we have 
\[  e(\mathrm{T} \mathcal{F}|S , \mathfrak{s}_0|\partial S) \cap [S] = \langle e(\mathrm{T}\mathcal{F} , \mathfrak{s}_0) , [S] \rangle. \]
Therefore, it is enough to show that
\[ \text{Ind}(\mathcal{F}, S) := e(\mathrm{T} \mathcal{F}|S , r|\partial S) \cap [S] = e(\mathrm{T} \mathcal{F}|S , \mathfrak{s}_0|\partial S) \cap [S].   \]
The relative Euler class $e(\mathrm{T} \mathcal{F}|S , \mathfrak{s}_0|\partial S)$ only depends on the homotopy class of $\mathfrak{s}_0|\partial S$ through non-vanishing sections. There is a homotopy from $\mathfrak{s}_0|\partial S$ to $r|\partial S$ through non-vanishing sections of $\mathrm{T} \mathcal{F}|\partial S$. This completes the proof. 
	
\end{proof}

\begin{remark}
	Let $M$ be a sutured 3-manifold, and $\mathcal{F}$ be a foliation of $M$ that is compatible with the sutured structure. Assume that some boundary component $T$ of $ M$ is a leaf of $\mathcal{F}$ and that $\chi(T)<0$. Then Definition \ref{def:ind} of the index sum $\text{Ind}(\mathcal{F}, \cdot)$ does not determine a well-defined cohomology class in $H^2(M , \partial M) \cong \text{Hom}(H_2(M , \partial M) , \mathbb{R})$. To see this, note that $T$ is a leaf of $\mathcal{F}$, and hence $|\text{Ind}(\mathcal{F} , T)| = |\chi(T)| > 0$, but $T$ represents the trivial homology class in $H_2(M , \partial M)$.
	\label{higher-genus-boundary}
\end{remark}

\subsection{Roussarie--Thurston general position, and Thurston's inequality}

Let $M$ be a compact oriented 3-manifold, and $\mathcal{F}$ be a transversely oriented foliation on $M$ such that each component of $\partial M$ is either a leaf of $\mathcal{F}$ or is transverse to $\mathcal{F}$. Let $S$ be a compact, properly embedded, oriented surface $S$ in $M$ such that each component of $\partial S$ is either contained in a leaf or is transverse to $\mathcal{F}$. Then $S$ can be isotoped to be in general position with respect to the foliation such that $S$ and $\mathcal{F}$ are transverse except at a finite number of saddle or center tangencies. 

If furthermore the foliation is taut, and the surface $S$ is incompressible and $\partial$-incompressible, then Roussarie \cite{roussarie1974plongements} and Thurston \cite{thurston1972foliations,thurston1986norm} showed that the surface $S$ can be isotoped such that each component of $S$ either is a leaf of $\mathcal{F}$ or is transverse to $\mathcal{F}$ except at finitely many points of saddle tangencies. Roussarie proved the statement for incompressible tori, and Thurston generalized it to any embedded incompressible surface. Gabai generalized it to the case of immersed incompressible surfaces and without any orientability assumption on the manifold and the foliation \cite{gabai2000combinatorial}. 

More generally, the Roussarie--Thurston general position holds when $(M,\gamma) $ is a sutured manifold. In this case, we assume that each component of $\partial S$ either 
\begin{enumerate}[i)]
	\item  is a leaf of $\mathcal{F}|\gamma$, or
	\item is transverse to $\mathcal{F}|\gamma$ and intersects $R(\gamma)$ in essential arcs, or
	\item intersects $R(\gamma)$ in an essential simple closed curve.

\end{enumerate}

Thurston's inequality generalizes to taut foliations on sutured manifolds and admissible surfaces as
\[ |\text{Ind}(\mathcal{F}, S)| \leq |\chi_s(S)|. \] 

\begin{remark}
	Both general positions that we mentioned hold for $C^{\infty,0}$ foliations as well \cite{solodov1984components,gabai2000combinatorial} although the original argument was for $C^2$ foliations.
\end{remark}

\subsection{Relative Euler class using the canonical trivialization of the tangent bundle of a torus} We can use the following construction to define a relative Euler class \cite[Section 3]{juhasz2008floer}. Let $(M,\gamma)$ be a \emph{balanced} sutured 3-manifold, meaning that $\chi(R_+(\gamma)) = \chi(R_-(\gamma))$, $M$ has no closed component, and the map $\pi_0(A(\gamma)) \rightarrow \pi_0(\partial M)$ is surjective. Let $\mathcal{F}$ be a foliation of $M$ compatible with the sutured structure. Pick a Riemannian metric on $M$. Denote by $v_0$ a nowhere vanishing vector field along $\partial M$ that points outward along $R_+(\gamma)$, points inward along $R_-(\gamma)$, and on $\gamma$ is the gradient of the height function $s(\gamma) \times I \rightarrow I$. The restriction of the tangent bundle of $\mathcal{F}$ to $\partial M$ can be identified with the orthogonal plane bundle $v_0^\perp$ to $v_0$. It is easy to see that the plane bundle $v_0^\perp$ is trivial (or equivalently has a section) if and only if for every component $F$ of $\partial M$ the equality $\chi(F \cap R_+(\gamma)) = \chi(F \cap R_-(\gamma))$ holds; such sutured manifolds are called \emph{strongly balanced} by Juh\'{a}sz. See Juh\'{a}sz \cite[Proposition 3.4]{juhasz2008floer}. \\

Assume that $(M, \gamma)$ is strongly balanced. Let $t$ be a trivialization of $v_0^\perp$. Then it makes sense to talk about the relative Euler class $e(\mathrm{T}\mathcal{F}, t) \in H^2(M , \partial M ; \mathbb{Z})$, which is defined as $e(\mathrm{T} \mathcal{F}, \mathfrak{s}_0)$ where $\mathfrak{s}_0$ is any nonzero constant section of $v_0^\perp$ under the identification $t \colon v_0^\perp \xrightarrow[]{\cong} \partial M \times \mathbb{R}^2$. Let $b$ be a component of $\partial S$ with the induced orientation from $S$. If $r$ is the vector field as in Definition \ref{def:ind}, then define $\text{rot}(b, t)$ to be the rotation number of $r|b$ with respect to $t$. Define $\text{rot}(S , t)$ as sum of $\text{rot}(b, t)$ where $b$ varies over the components of $\partial S$. The relative Euler class $e(\mathrm{T}\mathcal{F}, t)$ depends on the choice $t$ of trivialization, and its evaluation on an admissible surface $S$ is equal to the index sum in Equation (\ref{index-sum}) minus the correction term $\text{rot}(S,  t)$

\begin{eqnarray}
\langle e(\mathrm{T}\mathcal{F}, t), [S]\rangle = \Big( \sum_{\substack{\Text{tangent }\\ \Text{points }p}} \text{sign}(p) \cdot i(p) \Big) - \text{rot}(S , t).
\end{eqnarray}

This is because, by the naturality of the Euler class
\[ \langle e(\mathrm{T} \mathcal{F} , t) , [S] \rangle = \langle e(\mathrm{T}\mathcal{F}|S , t |\partial S) , [S] \rangle,  \]
and moreover 
\[ \langle e(\mathrm{T}\mathcal{F}|S , r |\partial S) , [S] \rangle - \langle e(\mathrm{T}\mathcal{F}|S , t |\partial S) , [S] \rangle = \text{rot}(S,t). \]

\begin{remark}
 Let $\nu_S$ be the positive unit normal vector to $S$, and assume that $\nu_S$ is nowhere parallel to $v_0$ along $\partial S$. This holds for generic $S$. Let $p(\nu_S)$ be the projection of $\nu_S$ onto $v_0^\perp$. Then the restriction of $p(\nu_S)$ to $\partial S$ is non-vanishing. In \cite{juhasz2008floer}, Juh\'{a}sz defines $\text{rot}(b , t)$ for a component $b$ of $\partial S$ as the rotation number of $p(\nu_S)|b$ with respect to the trivialization $t$. The above definition is equivalent to that of Juh\'{a}sz. This is because $(p(\nu_S(x)), r(x))$ is a positive basis for $T_x(\mathcal{F})$ for each $x \in b$; see Figure \ref{normal-vectors}, left. 
\end{remark}

The tangent bundle of an oriented two-dimensional torus, $T$, has a canonical trivialization up to homotopy, obtained as follows. Split $T$ as a product of two oriented circles. Trivialize the tangent bundles of oriented circles according to their orientations, and equip the torus $T$ with the product trivialization. Up to homotopy, this trivialization is invariant under Dehn twists along $S^1 \times \text{point}$ and $\text{point} \times S^1$, and hence does not depend on the splitting $T = S^1 \times S^1$. See Turaev \cite[Page 163]{turaev2012torsions}. Note we do not claim that there is a unique trivialization of the tangent bundle up to homotopy; we merely say that one homotopy class of trivializations is distinguished from the others. Likewise, the tangent bundle of an oriented annulus $A$ has a canonical trivialization coming from a splitting $A = I \times S^1$.

There is an intuitive characterization of the canonical trivialization of the tangent bundle of an oriented torus $T$ as follows: Given a trivialization $t$, there is a homomorphism 
\[ \phi_t \colon \pi_1(T) \longrightarrow \mathbb{Z} \]
such that for every smooth based loop $\gamma$, $\phi_t([\gamma])$ is the rotation number of the unit tangent vector $\dot{\gamma}$ with respect to $t$. For the canonical trivialization $t$, the homomorphism $\phi_t$ is trivial.\\

\begin{figure}
\labellist
\pinlabel $R_+(\gamma)$ at 40 67
\pinlabel $R_-(\gamma)$ at 105 67
\endlabellist	
	
\centering
\includegraphics[width=1.5 in]{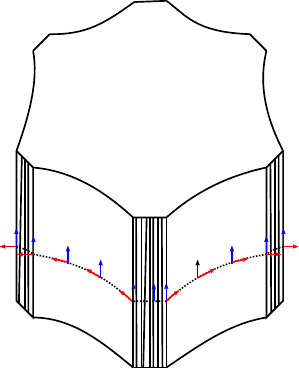}
\caption{The canonical (up to homotopy) trivialization of the restriction of the tangent bundle of the foliation to a torus boundary component containing sutures. In the above, $v_1$ is the vertical vector, and $v_2$ points inside the manifold along the middle annulus suture.}
\label{trivialization}
\end{figure}

\begin{definition}
Let $(M, \gamma)$ be a sutured 3-manifold with toroidal boundary, and $t_0$ be the trivialization of $\mathrm{T}\mathcal{F}|\partial M$, canonical up to homotopy, defined as follows. Pick a Riemannian metric on $M$.
\begin{enumerate}[i)]
	\item Along any torus component of $\partial_\pitchfork M$, $t_0$ is the trivialization associated to the outward section of $\mathrm{T} \mathcal{F}$. In other words, the trivialization $t_0$ picks the basis $\{o, i(o) \}$ where $o$ is the outward pointing section and $i(o)$ is obtained by rotating $o$ by angle $\frac{\pi}{2}$ inside $\mathrm{T} \mathcal{F}$.
	\item Along any torus component of $\partial_\tau M$ with the induced orientation from $\mathcal{F}$, let $t_0$ be the canonical trivialization of the tangent bundle of the oriented torus as described previously.
	\item Let $\mathbb{T}$ be a torus component of $\partial M$ with $2k>0$ sutures. Denote the tangential annuli on $\mathbb{T}$ by $A_1, A_2, \cdots, A_{2k}$ and equip them with the induced orientations from $\mathcal{F}$. Pick a simple closed curve $\alpha \subset \mathbb{T}$ parallel to a suture, and give $\alpha$ an orientation. Let $v_1$ be the unit tangent vector to $\mathcal{F}|(\mathbb{T} \cap \partial_{\pitchfork}M)$ whose orientation is coherent with that of $\alpha$. Split each oriented $A_j$ as a product $S_j^1 \times I_j$ where $S^1_j$ (respectively $I_j$) is an oriented circle (respectively interval) such that the orientation of $S_j^1$ is coherent with that of $\alpha$. On each $A_j$ pick the trivialization $\{ v_1 , v_2\}$ corresponding to its splitting $S_j^1 \times I_j$. At this point $v_1$ is defined on all of $\mathbb{T}$, and so this trivialization can be extended to a trivialization $\{ v_1, v_2 \}$ of $\mathrm{T}\mathcal{F}|\mathbb{T}$. The restriction of $v_2$ to each component of $\mathbb{T} \cap \partial_{\pitchfork} M$ will be either the inward pointing or the outward pointing section of $\mathrm{T} \mathcal{F}$. See Figure \ref{trivialization}.
\end{enumerate}
We call this trivialization \emph{canonical}, which is well-defined up to homotopy.
\end{definition}

\begin{definition}
Let $(M,\gamma)$ be a sutured 3-manifold with toroidal boundary, and $\mathcal{F}$ be a foliation of $M$ that is compatible with the sutured structure. Define the relative Euler class as $e(\mathcal{F}, t_0) \in H^2(M , \partial M; \mathbb{Z})$, where $t_0$ is the canonical trivialization of $\mathrm{T}\mathcal{F}|\partial M$.
\label{def:relative-euler-class}
\end{definition}  

\begin{prop}
Let $(M,\gamma)$ be a sutured 3-manifold with toroidal boundary, and $\mathcal{F}$ be a foliation of $M$ that is compatible with the sutured structure of $M$. Then $\text{Ind}(\mathcal{F}, \cdot)$ defines a cohomology class in $H^2(M, \partial M)$ which coincides with the real Euler class relative to the canonical trivialization of $\mathrm{T}\mathcal{F}|\partial M$.
\end{prop}

\begin{proof}
For every admissible surface $S$ and component $b$ of $\partial S$, the rotation number $\text{rot}(b,  t_0)$ is equal to zero.
\end{proof}

\subsection{Parity condition}
The following parity condition goes back at least to Wood \cite{wood1969foliations}.

\begin{prop}[Parity condition]
Let $M$ be a closed orientable 3-manifold, and $\xi$ be a transversely oriented plane field on $M$. Then the Euler class of $\xi$ lies in $2 H^2(M; \mathbb{Z})$.
\label{parity}
\end{prop}

\begin{proof}
Since $\xi$ is transversely oriented, $\xi \oplus \epsilon \cong \mathrm{T}M$ where $\epsilon$ is the trivial line bundle. Therefore the Stiefel--Whitney classes $w_2(\xi)$ and $w_2(\xi \oplus \epsilon) = w_2(M)$ are equal. As $M$ is parallelizable, we have $w_2(M)=0$ implying that $w_2(\xi)=0$. But $w_2(\xi)$ is the mod $2$ reduction of $e(\xi)$, hence $e(\xi) \in 2H^2(M; \mathbb{Z})$. For an alternative proof, see the proof of Theorem \ref{Wood}.
\end{proof}

\subsection{Standard taut foliations by stack of saddles}

\begin{figure}
	\centering
	\begin{minipage}{4.2 cm}
		\centering
		\includegraphics[width=2 in]{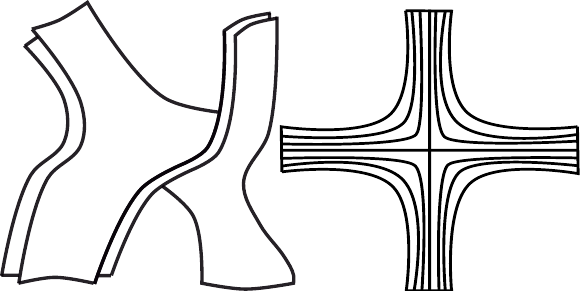}
	\end{minipage}
	\begin{minipage}{4.2 cm}
		\centering
		\includegraphics[width=1 in]{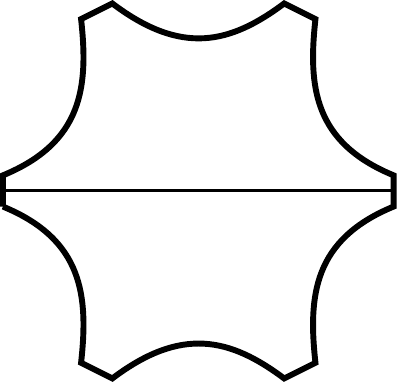}
	\end{minipage}
	\caption{Left: a foliation of $M$ by a stack of chairs. Middle: the induced foliation on a meridional disk of $M$. Right: a meridional disk of $\hat{N}$ is divided into two disks.}
	\label{stackofchairs}
\end{figure}

\begin{ex}
	Let $M$ be the sutured solid torus as in Example \ref{Ex-stack-of-chairs}. There is a taut foliation of $M$ by a stack of chairs obtained in the following way. Take an infinite stack of chairs and glue the top to the bottom by $180^\circ$ rotation. See Figure \ref{stackofchairs}, left. The reader who would like to see one possible equation defining the leaves of this foliation should consult \cite[Page 361]{candel2003foliations}. Recall that the meridional disk $F$ of $M$ has sutured Euler characteristic equal to $-1$. By the index sum formula, the relative Euler class of this foliation assigns $\pm 1$ to the meridional disk of $M$. Both $\pm 1$ can be realized by choosing appropriate transverse orientations. 
	
	One can define a similar foliation on the sutured solid torus $N$ with two sutures each of which goes three times around a longitude and once around the meridian. This foliation on $N$ looks like a stack of monkey saddles. In this case, the sutured Euler characteristic of a meridional disk is equal to $-2$, and the relative Euler class of the foliation assigns $\pm 2$ to the meridional disk. A similar foliation can be constructed for any sutured solid torus where the sutures are parallel essential simple closed curves on the boundary torus and intersect the meridian non-trivially, and it is called the \emph{standard taut foliation by a stack of generalized saddles}. Note that in these examples, the holonomy of each component of the transverse boundary is a \emph{shift map}, i.e. it has no fixed points except the interval endpoints.  
	\label{stackexample}
\end{ex}

\begin{ex}
	Let $\hat{N}$ be a sutured solid torus with six sutures, each of which is isotopic to a longitude. If $\hat{D}$ is the meridional disk of $\hat{N}$ then $\hat{N} \setminus \setminus \hat{D}$ is homeomorphic to $\hat{D} \times I$ where $I $ is an interval. Moreover, the image of the annuli sutures of $\hat{N}$ in $\hat{N} \setminus \setminus \hat{D}$ consists of $6$ disjoint parallel vertical rectangles $J_i \times I$ where $J_i$ are disjoint intervals in $\partial \hat{D}$. In particular, $\hat{N} \setminus \setminus \hat{D}$ and $N \setminus \setminus D$ are identical. See Figure \ref{solidtorus}, right. 
	
	We may construct a taut foliation on $\hat{N}$ whose relative Euler class assigns zero to $\hat{D}$ as follows. Divide $\hat{D}$ into two disks $\hat{D}_1$ and $\hat{D}_2$ each of which intersects the annuli sutures in $4$ arcs as in Figure \ref{stackofchairs}, right. This divides $\hat{N}$ into two solid tori $\hat{N}_1$ and $\hat{N}_2$ with meridional disks $\hat{D}_1$ and $\hat{D}_2$ respectively. Pick a taut foliation $\mathcal{F}_1$ on $\hat{N}_1$ (respectively $\mathcal{F}_2$ on $\hat{N}_2$) whose relative Euler class assigns $+1$ (respectively $-1$) to $\hat{D}_1$ (respectively $\hat{D}_2$) as constructed in Example \ref{stackexample}. Pick a positive (respectively negative) annulus component $T_1$ (respectively $T_2$) of the tangential boundary of $\hat{N}_1$ (respectively $\hat{N}_2$). Glue $\mathcal{F}_1$ to $\mathcal{F}_2$ by identifying $T_1$ with $T_2$ appropriately to obtain a transversely oriented foliation $\mathcal{F}$ of $\hat{N}$. It is easy to see that $\mathcal{F}$ is taut, and its relative Euler class assigns zero to $\hat{D}$ by the index sum formula. 
	\label{Ex-no-twist}
\end{ex}

\section{Construction of the 3-manifolds}  
In this section, we explain the construction of our counterexamples. Each of the constructed manifolds is obtained by a Dehn surgery on a fibered 3-manifold, but we should put some constraints on the monodromy of the fibered manifold and specify the surgery curve. \\

\textbf{Step I}: \textit{Definition of the fiber $S$ and curves $\alpha , \beta, \gamma \subset S$.}

\begin{enumerate}
	\item Let $S$ be a closed, orientable surface of genus $g \geq 3$.
	\item Let $\gamma$ be a non-separating simple closed curve on $S$, and $A$ be an annulus neighborhood of $\gamma$ in $S$ with $\partial A = \gamma _{+} \cup \gamma _{-}$. Pick an orientation on $\gamma$, which induces orientations on $\gamma _+$ and $\gamma _-$.
	\item Choose points $p_1, p_2, p_3 \subset \gamma_+$ that are cyclically ordered in the same direction as $\gamma$. Similarly choose points $q_1, q_2, q_3$ on $\gamma _-$.
    \item Choose disjoint oriented arcs $l_1, l_2, l_3 \subset A $ where $l_i$ connects $q_i $ to $p_i $. Likewise, choose disjoint oriented arcs $m_1, m_2, m_3 \subset A $ where $m_i$ connects  $q_{i+1} $ to $p_i $; the arithmetic on the indices is modulo $3$ throughout. Note the shift in indices. We choose them in such a way that $m_i$ is obtained by adding the oriented arc $q_{i+1}q_i$ to $l_i$ and then perturbing them to be disjoint and properly embedded in $A $. See Figure \ref{cylinder}.
    \begin{figure}
    	
    	\labellist
    	\small\hair 2pt
    	\pinlabel $p_1$ at 14 25
    	\pinlabel $p_2$ at 14 62
    	\pinlabel $p_3$ at 14 110
    	
    	\pinlabel $q_1$ at 164 24
    	\pinlabel $q_2$ at 164 63
    	\pinlabel $q_3$ at 164 110
    	
    	\pinlabel $p_1$ at 262 25
    	\pinlabel $p_2$ at 262 62	
    	\pinlabel $p_3$ at 262 110
    	
    	\pinlabel $q_1$ at 410 25
    	\pinlabel $q_2$ at 410 62
    	\pinlabel $q_3$ at 410 110
    	\endlabellist
    	
    	\centering
    	\includegraphics[width=4 in]{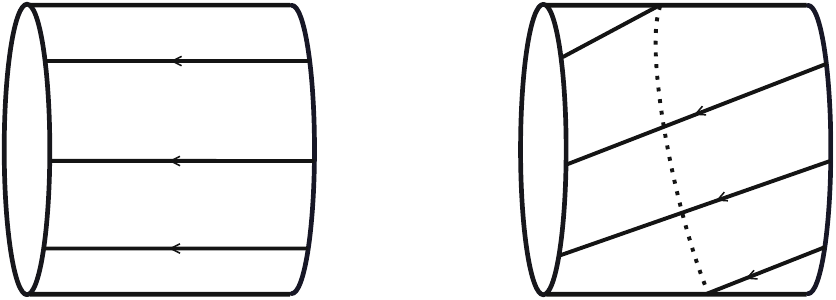}
    	\caption{Left: arcs $l_i \subset A$. Right: arcs $m_i \subset A$}
    	\label{cylinder}
    \end{figure}
    \item Let $\delta _1, \delta_2, \delta_3$ be disjoint oriented simple arcs in $S - A^\circ$ such that $\delta_i$ connects $p_i$ to $q_{i-1}$. Again note the shift in indices.
    \item Define $\alpha$ as the union of the six arcs $l_i$ and $\delta _{i}$ for $1 \leq i \leq 3$. The curve $\beta$ is defined as the union of $m_i$ and $\delta _{i} $ for $1 \leq i \leq 3$. 
\end{enumerate}    

  Note that the appropriate shifts in indices for $m_i$ and $\delta_i$ imply that $\alpha$ and $\beta$ are connected simple closed curves. Moreover, $\beta = \alpha - \gamma $ as oriented sum (oriented cut and paste). In particular, $\alpha$ and $\beta$ are not homologous to each other in $S$ since $[\gamma] \neq 0 \in H_1(S)$. \\

\textbf{Step II}: \textit{Definition of a homeomorphism $f$ of $S$, and its mapping torus $M_f$.}

\begin{lem}
	Let $S$ be a closed, orientable surface of genus $g=3$ or $g \geq 6$. For suitable choices of $\gamma$, $l_i$, $m_i$ and $\delta_i$ as in Step I, there exists a homeomorphism $f$ of $S$ such that
	\begin{enumerate}
		\item $f$ is pseudo-Anosov;
		\item $f(\alpha) = \beta$;
		\item $\Text{rank}(H_2(M_f))=1$, where $M_f$ is the mapping torus of $f$.
	\end{enumerate}
	\label{pseudoanosov}
\end{lem}

\begin{proof}
	The condition $\Text{rank}(H_2(M_f))=1$ is equivalent to requiring that the map 
	\[  (f_* - Id) \colon H_1(S) \longrightarrow H_1(S)\] 
	has trivial kernel \cite[Example 2.48]{hatcher2002algebraic}. We first construct the map for a surface of genus $3$ and from there it is clear how to generalize it for arbitrary genera $g \geq 6$. Let $\gamma$, $\alpha$ and $\beta$ be the oriented curves shown in Figures \ref{arcscurves-alpha} and \ref{arcscurves-beta}. The curves $\gamma$, $\alpha$ and $\beta$ are redrawn in Figure \ref{alphabeta} to make them more visible; the arcs $\delta_i$ can be thought of as the three pieces of $\alpha - \alpha \cap \beta$. It is easy to see that $\alpha$, $\beta$ and $\gamma$ are non-separating. We use Penner's construction of pseudo-Anosov maps to define $f$. For any simple closed curve $\eta$, let $\tau_{\eta}$ be the positive Dehn twist around $\eta$. Our convention is the right handed twist.

	Let $a_1 , \cdots , a_4$ (negative twists) and $b_1, \cdots , b_3$ (positive twists) be the filling system of curves shown in the left side of Figure \ref{pennercurves}. Define the maps $f$ as
	\[ f := \tau_{b_2} \circ \tau_{a_2}^{-1} \circ \tau_{a_3}^{-1} \circ \tau_{b_1} \circ \tau_{b_3} \circ \tau_{a_1}^{-1} \circ \tau_{a_4}^{-1}. \] 
	By Penner's construction, the map $f$ is pseudo-Anosov. Since the curves $a_1 , a_4 , b_1,b_3$ are disjoint from $\alpha$, we have
	\[ f(\alpha) =  \tau_{b_2} \circ \tau_{a_2}^{-1} \circ \tau_{a_3}^{-1}(\alpha) = \beta. \]
	
%
%
%

\begin{figure}
	\labellist
	\pinlabel $q_1$ at 78 105
	\pinlabel $p_1$ at 100 114
	\pinlabel $q_3$ at 135 62
	\pinlabel $p_3$ at 116 44
	\pinlabel $q_2$ at 75 52
	\pinlabel $p_2$ at 79 72
	\pinlabel $\delta_3$ at 67 0
	\pinlabel $l_1$ at 53 122
	\pinlabel $\delta_3$ at 408 15
	\pinlabel $\alpha$ at 117 97
	\pinlabel $\gamma_{+}$ at 115 123
	\pinlabel $\gamma_{-}$ at 115 145
	\endlabellist
	
	\centering
	\includegraphics[width= 4 in]{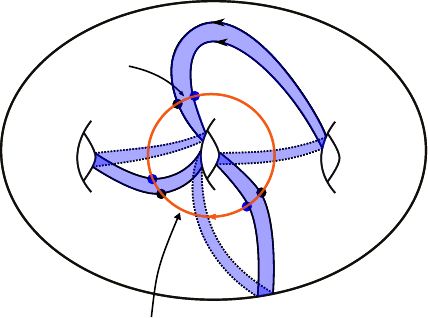}
		\caption{The relative position of the curve $\alpha$ with respect to the curves $\gamma_{+}, \gamma_{-}$, the points $p_i$, $q_i$, and the arcs $m_i$, $l_i$, and $\delta_i$ is shown. The shaded region is the annulus $A$ cobounded by $\gamma_{+}$ and $\gamma_{-}$.}
	\label{arcscurves-alpha}
\end{figure}

\begin{figure}
	\labellist
	\pinlabel $q_1$ at 78 95
	\pinlabel $p_1$ at 100 104
	\pinlabel $q_3$ at 132 50
	\pinlabel $p_3$ at 117 38
	\pinlabel $q_2$ at 75 46
	\pinlabel $p_2$ at 78 63
	\pinlabel $\beta$ at 127 91
	\pinlabel $\gamma_{+}$ at 115 115
	\pinlabel $\gamma_{-}$ at 115 135
	\pinlabel $m_1$ at 23 5
	\pinlabel $\delta_3$ at 167 3
	\endlabellist
	
	\centering
	\includegraphics[width= 4 in]{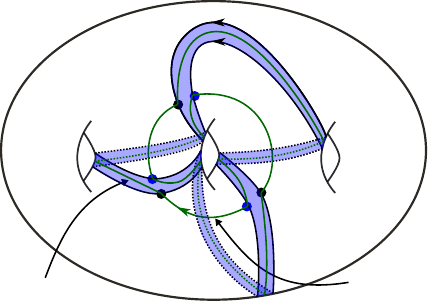}
	\caption{The curve $\beta$}
	\label{arcscurves-beta}
\end{figure}

	\begin{figure}
	\centering
	\labellist
	\pinlabel $\gamma$ at 315 45
	\pinlabel $\alpha$ at 100 45 
	\pinlabel $\beta$ at 57 30
	\endlabellist
	
	\includegraphics[width=5 in]{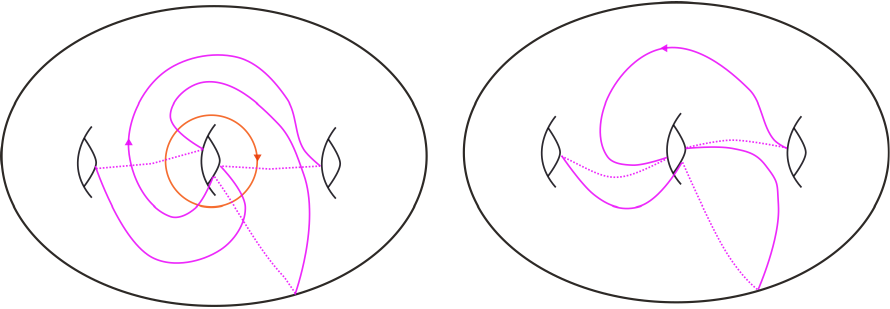}
	\caption{The curves $\alpha$, $\beta$ and $\gamma$ on $S$}
	\label{alphabeta}
\end{figure}

	\begin{figure}
		\labellist
		\pinlabel $a_1$ at 190 55
		\pinlabel $b_1$ at 160 40
		\pinlabel $a_2$ at 130 60
		\pinlabel $b_2$ at 100 100
		\pinlabel $a_3$ at 70 60 
		\pinlabel $b_3$ at 40 40
		\pinlabel $a_4$ at 10 60
		
		\pinlabel $r_1$ at 415 60
		\pinlabel $s_1$ at 385 45 
		\pinlabel $r_2$ at 355 60 
		\pinlabel $s_2$ at 325 45
		\pinlabel $r_3$ at 295 60
		\pinlabel $s_3$ at 265 45
		\endlabellist

		\centering
		\includegraphics[width=5 in]{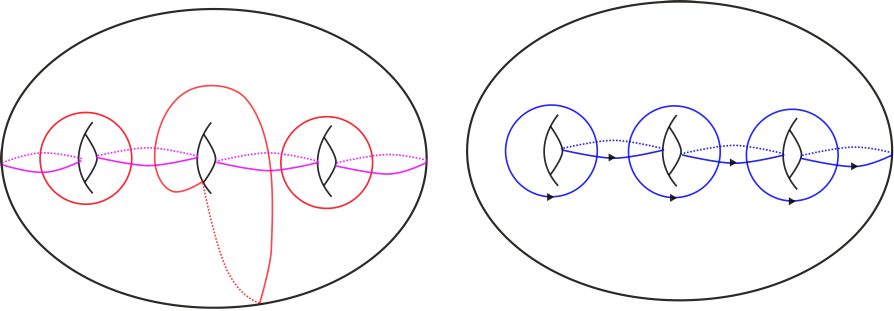}
		\caption{Left: Penner curves. Right: a basis for the first homology}
		\label{pennercurves}
	\end{figure}

	We would like to show that $\det(f_* -Id)$ is nonzero. Choose a basis $r_1, s_1 , \cdots , r_3, s_3$ for $H_1(S)$ as in the right side of Figure \ref{pennercurves}. One can directly see that the action of $f$ on homology is represented by the matrix
	\begin{eqnarray*}
		W = \begin{bmatrix}
			0  & 1  & 2 & -1  & -2 &1 \\
			-1  & 2  & 1 & 0 & 0  & 0\\
			-2  & 4  & 4 & -2 & -2 & 1\\
			1  & -2  & -2 & 2 & 2 & -2 \\
			0  &0    & 0 & 1 & 2 & -3 \\
			0 & 0    & 0 & 0 & -1 & 2 \\
		\end{bmatrix},
	\end{eqnarray*}
and that $\det(W-Id) \neq 0$. This finishes the proof for the genus $3$ surface. For $g \geq6$, add extra handles to the left side of the picture and add suitable curves to complete the previous system of filling curves as in Figure \ref{pennercurvesgeneral} and Figure \ref{homologygeneral}. 
	\begin{figure}
		\labellist
		\pinlabel $a_1$ at 385 52 
		\pinlabel $b_1$ at 355 35
		\pinlabel $a_2$ at 315 55 
		\pinlabel $b_2$ at 290 90 
		\pinlabel $a_3$ at 263 55 
		\pinlabel $b_3$ at 240 35
		\pinlabel $a_4$ at 208 55
		\pinlabel $b_4$ at 183 35
		\pinlabel $a_5$ at 150 55
		\pinlabel $b_5$ at 125 35
		\pinlabel $a_{g+1}$ at 17 55 
		\pinlabel $b_g$ at 50 35 
		\endlabellist
		
		\centering
		\includegraphics[width=4 in]{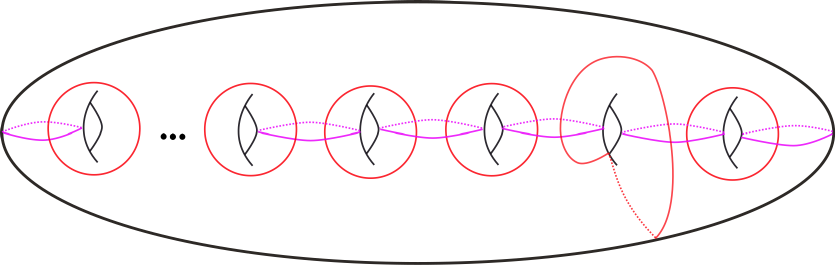}
		\caption{Penner curves for the general case }
		\label{pennercurvesgeneral}
	\end{figure}
	
	\begin{figure}
		\centering
		\labellist
		\pinlabel $r_1$ at 385 50
		\pinlabel $r_2$ at 325 50 
		\pinlabel $r_3$ at 265 50 
		\pinlabel $r_4$ at 205 50 
		\pinlabel $r_5$ at 147 50 
		
		\pinlabel $s_1$ at 355 35 
		\pinlabel $s_2$ at 295 35 
		\pinlabel $s_3$ at 235 35 
		\pinlabel $s_4$ at 175 35 
		\pinlabel $s_5$ at 115 35 
		
		\pinlabel $s_g$ at 50 35
		\endlabellist
		
		\includegraphics[width=4 in]{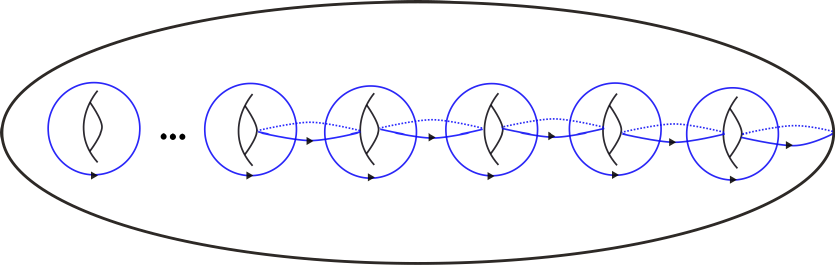}
		\caption{A basis for the first homology }
		\label{homologygeneral}
	\end{figure}
	
	The map $f$ is defined similarly. First we do negative twists around $\{ a_1 , \cdots , a_{g+1}\} \setminus \{a_2 , a_3\}$. Then we do positive twists around $\{b_1, \cdots , b_g \} \setminus \{ b_2\}$, and at the end we do negative twists around $a_2 , a_3$ followed by a positive twist around $b_2$. Again computation shows that for $g \geq 6$ the action of $f$ on homology is represented by the following matrix $V$. Here the empty entries are zero and $*$ shows a repeating pattern. The patterns correspond to 
	\begin{align*}
	 r_i & \longmapsto s_{i-1}+r_i -s_i, &5 \leq i \leq g-1, \\
	 s_i & \longmapsto -s_{i-1} -r_i +3s_i +r_{i+1}-s_{i+1}, &5 \leq i \leq g-1,
	 \end{align*}
	as well as
	\[ s_g \longmapsto -s_{g-1} - r_g + 3 s_g - (r_1+ r_2 + \cdots + r_g).  \]
	
	\[
	V = \left[\begin{array}{cccccccc|cc|cc|cc|cc}
	0 & 1 & 2 & -1 & -2 & 2 & 1 & -1 & & & &&&&0 &-1 \\
	-1 & 2 & 1 & 0 & 0 &0 & 0 & 0 & &&&&&& 0 & 0 \\ 
	-2 & 4 & 4 & -2 & -2 & 2 & 1 & -1 & &&&&&& 0 & -1 \\ 
	1 & -2 & -2 & 2 & 2 & -2 & -1 & 1 &&&&&&& 0 & 0 \\
	&&&  1 & 2 & -2 & -1 & 1 & &&&&&& 0 & -1 \\
	&&&& -1 & 2 & 1 & -1 &&&&&&& 0 & 0 \\
	&&&&& 1 & 1 & -1 &&&&&&& 0 & -1 \\
	&&&&& -1 & -1 & 3 & 1 & -1 &&&&& 0 & 0 \\ \hline
	&&&&&&& 1 & 1 & -1 &&&&& * & * \\
	&&&&&&& -1 & -1 & 3 &* & * &&&* & * \\ \hline
	&&&&&&&&& * & * & * &&& 0 & -1 \\
	&&&&&&&&&  *&  *&* &1&-1& 0 & 0 \\ \hline 
	&&&&&&&&&&&1 &1 & -1 & 0 & -1 \\
	&&&&&&&&&&&-1&-1&3&1&-1 \\ \hline 
	&&&&&&&&&&&&&1 & 1 & -2 \\
	&&&&&&&&&&&&&-1 & -1 & 3 \\  
	\end{array}\right].
	\]
	
	By putting the matrix $V - Id$ in the row reduced form, one can see that it is invertible. We do this in detail. Put the first eight rows in the row reduced form. After this, the first eight rows of $V-Id$ change to
	\[ 
	\left[\begin{array}{cccccccc|cc|cc|cc|cc}
	1 &  &  &  &  &  &  &  &-1 &1 & &&&&0 &3 \\
	& 1 & & & & &  &  & &&&&&& 0 & 1 \\ 
	&  & 1 & & & & &  &-1 &1&&&&& 0 & 2 \\ 
	&  &  &  1&  &  &  &  &&&&&&& 0 & 2 \\
	&&&  & 1 &  &  &  & -1&1&&&&& 0 & 4 \\
	&&&&  & 1 &  &  &&&&&&& 0 & 3 \\
	&&&&&  & 1 &  &-1&1&&&&& 0 & 5 \\
	&&&&&  &  & 1 &  &  &&&&& 0 & 4 \\ \hline
	\end{array}\right].
	\]
	Let us denote the $i$-th row of the matrix by $R_i$. For $5 \leq i \leq g$, do the following moves
	\begin{enumerate}
		\item Replace $R_{2i-1}$ with $R_{2i-1} - R_{2i-2}$.
		\item Replace $R_{2i}$ with $R_{2i} + R_{2i-2}$.
		\item Switch $R_{2i}$ and $R_{2i-1}$.
		\item Replace $R_{2i}$ with $- R_{2i}$.
	\end{enumerate}

	This process makes the matrix upper triangular, with all the entries on the diagonal equal to $\pm 1$ except the last one that is equal to $-(g+1)$. Therefore, for $g \geq 6$ the determinant of $V-Id$ is equal to $g+1$ up to sign. This shows that $V-Id$ has trivial kernel.
	
\end{proof}

\begin{remark}
	We avoided the case $g=4,5$ since it required a separate computation and it was not necessary for the proof of the main theorem, otherwise one can do similar constructions.
\end{remark}

\textbf{Step III}: \textit{Definition of the manifold $M = M(S ,\gamma, \ell, m,\delta, f)$}, where $\ell := \{ \ell_1 , \ell_2, \ell_3 \}$ and similarly for $m$ and $\delta$.\\

\begin{figure}
\labellist
\pinlabel $l_1$ at 83 135
\pinlabel $l_2$ at 133 120
\pinlabel $l_3$ at 183 105
\pinlabel $m_1$ at 96 65
\pinlabel $m_2$ at 143 50
\pinlabel $m_3$ at 250 25
\pinlabel $m_3$ at 0 105
\pinlabel $p_1$ at 72 30
\pinlabel $p_2$ at 120 18
\pinlabel $p_3$ at 168 6
\pinlabel $q_1$ at 110 105
\pinlabel $q_2$ at 160 91
\pinlabel $q_3$ at 211 79
\pinlabel {$A \times \frac{1}{4}$} at 280 60
\pinlabel {$A \times \frac{3}{4}$} at 280 90

\endlabellist
\centering
\includegraphics[width=4 in]{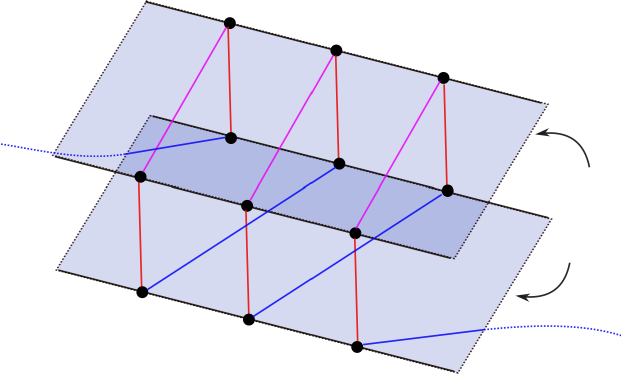}
\caption{The curve $\partial D$ is identified with the union of $l_i \times \{ \frac{3}{4} \} , m_i \times \{ \frac{1}{4} \} , p_i \times [\frac{1}{4}, \frac{3}{4}]$, and $q_i \times [\frac{1}{4}, \frac{3}{4}]$ (blue, red and purple curves) on the boundary of $M_f - U^\circ$. For simplifying the figure, $l_i \times \{ \frac{3}{4} \}$ is labelled with $l_i$; similarly for $m_i \times \{\frac{1}{4} \}$.} 
\label{boundaryofD}
\end{figure}

Let $S$ be a closed, orientable surface of genus $g=3$ or $g \geq 6$, and let $M_f$ be a fibered 3-manifold with fiber $S$ and monodromy $f \colon S \longrightarrow S$ where $f$ is as in \textbf{Step II}. Therefore, $M_f$ is obtained from $S \times [0,1]$ by identifying $(x,1)$ with $(f(x),0) $ for every $x \in S$. Recall that $A$ is an annulus neighborhood of $\gamma$ in $S$ with $\partial A = \gamma_+ \cup \gamma_-$. Let $U$ be the solid torus $A \times [\frac{1}{4} , \frac{3}{4}]$. The manifold $M$ is obtained from $M_f$ by removing the interior of $U$, denoted by $U^\circ$, and attaching the solid torus $N$ (Figure \ref{solidtorus})
\[ M = (M_f - U^\circ) \cup N, \] 
in the following way. Let $D$ be the meridional disk of $N$. The gluing of $N$ to $M_f - U^\circ$ is chosen such that $\partial D$ is identified with the union of $l_i \times \{ \frac{3}{4} \}$, $m_i \times \{ \frac{1}{4} \}$, $p_i \times [\frac{1}{4}, \frac{3}{4}]$, and $q_i \times [\frac{1}{4}, \frac{3}{4}]$ for $1 \leq i \leq 3$. See Figure \ref{boundaryofD}. Note $M_f-U^\circ$ can be considered as a cornered manifold with transverse boundary $\{\gamma_+ \cup \gamma_- \} \times (\frac{1}{4},\frac{3}{4})$, tangential boundary $A^\circ \times \{ \frac{1}{4}, \frac{3}{4} \}$, and concave corners along $\partial A \times \{ \frac{1}{4}, \frac{3}{4} \}$. Considering $N$ as a sutured manifold, the gluing map from $\partial N$ to $\partial (M_f - U^\circ)$ maps the transverse boundary to the transverse boundary, and the tangential boundary to the tangential boundary.

\label{construction-of-manifolds}
\section{Properties of the constructed manifolds}
\label{section:properties-of-manifolds}

In this section, we compute the Thurston norm of the constructed manifolds in Section \ref{construction-of-manifolds}, and show that the constructed manifolds are hyperbolic.
 
\begin{lem}[The Thurston norm of constructed manifolds]
Let $M$ be one of the manifolds constructed in Section \ref{construction-of-manifolds}. The second homology group $H_2(M)$ has rank two and the unit balls for the Thurston norm and the dual Thurston norm of $M$ are as in Figure \ref{unitballs}.
\label{Thurston-norm}
\end{lem} 
 
\begin{proof}\ \\
	\textit{Computation of the first Betti number}: Since $H_2(M_f)$ has rank $1$, an application of the Mayer--Vietoris sequence shows that $H_2(M)$ has rank at most $2$. More precisely, consider the following exact sequence for $M_f = (M_f-U^\circ) \cup U$
\begin{eqnarray*}
\dots \rightarrow H_1(\partial U) \rightarrow H_1(M_f-U^\circ)\oplus H_1(U) \rightarrow H_1(M_f) \rightarrow 0.
\end{eqnarray*}
This implies that
\begin{eqnarray*}
\Text{rank} \big(H_1(M_f-U^\circ) \big)+\Text{rank} \big( H_1(U) \big) \leq \Text{rank} \big( H_1(\partial U) \big) + \Text{rank} \big( H_1(M_f) \big)
\end{eqnarray*}
\begin{eqnarray*}
\Rightarrow \Text{rank} \big( H_1(M_f-U^\circ) \big) \leq 2.
\end{eqnarray*}
Note that $M_f - U^\circ$ and  $M - N^\circ$ are homeomorphic. Hence, we also have 
\begin{eqnarray}
\Text{rank} \big( H_1(M-N^\circ) \big) \leq 2.
\label{upper-bound}
\end{eqnarray}
Likewise, looking at the exact sequence for $M = (M-N^\circ) \cup N$ we have
\begin{eqnarray*}
\dots \rightarrow H_1(\partial N) \rightarrow H_1(M-N^\circ)\oplus H_1(N) \xrightarrow{\Phi} H_1(M) \rightarrow 0.
\end{eqnarray*}
This implies that
\begin{eqnarray}
H_1(M)\cong \frac{H_1(M-N^\circ)\oplus H_1(N)}{\ker (\Phi)},
\label{quotient}
\end{eqnarray}
where $\Phi$ is defined as $\Phi(x,y) = x+y$. We claim that 
\begin{eqnarray}
\Text{rank}(\ker(\Phi)) \geq 1.\label{lower-bound}
\end{eqnarray} 
To see this, let $y \in H_1(N)$ represent the core curve of $N$ and $x \in H_1(M-N^\circ)$ be the element that corresponds to $y$ via the attaching map between $N$ and $M-N^\circ$. Then 
\[  (0,0) \neq (-x,y) \in \ker (\Phi), \]
and no power of $(-x,y)$ is zero either. Hence $\Text{rank}(\ker (\Phi))$ is at least one. The isomorphism (\ref{quotient}) implies that
\begin{eqnarray*}
\Text{rank} \big( H_1(M) \big) = \Text{rank} \big( H_1(M-N^\circ) \big) + \Text{rank} \big( H_1(N) \big) - \Text{rank}\big(\ker (\Phi) \big) \leq 2 +1 -1 =2,
\end{eqnarray*}
where we have used Inequalities (\ref{upper-bound}) and (\ref{lower-bound}) for the last implication. This finishes the proof of the upper bound for the first Betti number of $M$.

On the other hand, the satisfied Conditions (1)--(3) by $f$ in the statement of Lemma \ref{pseudoanosov} imply that $H_2(M)$ has rank at least $2$. The surface $S \times \{ 1 \}$ is a nontrivial second homology class, as it admits a curve intersecting it transversely and exactly once. A second surface $F$ can be obtained as follows. Let $D$ be the meridional disk of $N$. We can assume that $\partial D$ is the union of $l_i \times \{ \frac{3}{4} \}$, $m_i \times \{ \frac{1}{4} \}$, $p_i \times [\frac{1}{4}, \frac{3}{4}]$, and $q_i \times [\frac{1}{4}, \frac{3}{4}]$ for $1 \leq i \leq 3$. See Figure \ref{boundaryofD}. Attach the three bands $\delta_i \times [\frac{1}{4},\frac{3}{4}]$ to $D$ to obtain a surface $F_0$ with two boundary components $\alpha \times \{ \frac{3}{4} \}$ and $-\beta \times \{ \frac{1}{4} \}$ and with Euler characteristic equal to $-2$, where $- \beta$ denotes $\beta$ with the opposite orientation. Since $f$ sends $\alpha$ to $\beta$, we can close the surface $F_0$ to get a closed orientable surface $F$ of genus two. More precisely, the surface $F$ is obtained from the union of $D$, the three mentioned bands, and the two vertical annuli $\alpha \times [\frac{3}{4},1]$, $\beta \times [0,\frac{1}{4}]$ by identifying $\alpha \times \{1 \}$ with $\beta \times \{ 0 \}$. Note that the homology classes $[S]$ and $[F]$ are linearly independent, since there is a curve that intersects $S$ (respectively $F$) transversely once and is disjoint from $F$ (respectively $S$). For example, the first curve corresponds to an appropriate vertical curve in $M_f$, and the second can be taken to be the core of $N$.\\

\textit{Computation of the Thurston norm}: We already know that $H_2(M)$ is generated by $S$ and $F$, where we have chosen fixed orientations on $S$ and $F$. These classes have Thurston norms at most $2g-2$ and $2$ respectively. Therefore 
\begin{eqnarray}
x([S]+[F]) \leq x([S])+x([F]) \leq -\chi(S)-\chi(F).
\label{thurston-norm-sum}
\end{eqnarray}
We exhibit taut foliations on $M$ and show that the unit balls for the Thurston norm and dual Thurston norm are as in Figure \ref{unitballs}. We will show that there is a taut foliation $\mathcal{F}_1$ of $M$ whose Euler class $e(\mathcal{F}_1)$ assigns the numbers $\chi(S)$ and $\chi(F)$ to the surfaces $S$ and $F$ respectively, implying that
\begin{eqnarray}
x([S])=-\chi(S),
\label{thurston-norm-one}
\end{eqnarray}
\begin{eqnarray}
x([F])=-\chi(F).  
\label{thurston-norm-two}
\end{eqnarray}
Hence
\[ \langle e(\mathcal{F}_1) , [S]+[F] \rangle = \langle e(\mathcal{F}_1) , [S] \rangle + \langle e(\mathcal{F}_1) , [F] \rangle = \chi(S)+\chi(F) \leq - x([S]+[F]), \]
where we used Inequality (\ref{thurston-norm-sum}) for the last implication. But we know that the dual norm of $e(\mathcal{F}_1)$ is at most one, therefore the equality should happen in the above, i.e.
\begin{eqnarray}
x([S]+[F]) = -\chi(S)-\chi(F).
\label{thurston-norm-three}
\end{eqnarray}  

To construct such a taut foliation $\mathcal{F}_1$ on $M$, decompose $M$ into three pieces: 
\[ N, \hspace{3mm} \big(S \times [\frac{3}{4},1] \big) \cup \big( S \times [0, \frac{1}{4}] \big), \hspace{3mm} (S-A^\circ) \times [\frac{1}{4}, \frac{3}{4}]. \]
On the first piece, choose a standard taut foliation of $N$ by a stack of monkey saddles whose relative Euler class assigns $\chi_s (D) =-2$ to the meridional disk of $N$. See Example \ref{stackexample}. Note that the holonomy of the transverse boundary of $N$ is a shift map. For the second piece, choose the product foliation. For the third piece, choose a foliation transverse to $[\frac{1}{4}, \frac{3}{4}]$ factor that has appropriate shift holonomies on $\gamma_+ \times [\frac{1}{4}, \frac{3}{4}]$ and $\gamma_- \times [\frac{1}{4}, \frac{3}{4}]$. Such a foliation exists by Lemma \ref{transversefoliations}. Glue these foliations along their common boundaries to obtain a foliation of $M$. See Figure \ref{product-norm}. The constructed foliation is taut since its compact leaves, essentially the surface $S$, have closed transversals. The Euler class of $\mathcal{F}_1$ assigns the numbers $\chi(S)$ and $\chi(F)$ to the surfaces $S$ and $F$ respectively. This is clearly true for $S$, since it is a leaf. To see that it also holds for $F$, recall that the number assigned to $F$ can be computed from an index sum formula. The surface $F$ can be obtained from $D$ by adding the three bands $\delta_i \times [\frac{1}{4}, \frac{3}{4}]$ as well as the two annuli $\alpha \times [\frac{3}{4},1]$ and $\beta \times [0, \frac{1}{4}]$, and then gluing the two boundary curves $\alpha \times \{1 \}$ and $\beta \times \{ 0 \}$ together. Since the induced foliations on the bands and annuli are product foliations by the construction, there is no tangency on the bands and annuli, and so they do not contribute to the index sum. Hence 
\[ \langle e(\mathcal{F}_1), F \rangle = \text{Ind}(\mathcal{F}_1, F) = \text{Ind}(\mathcal{F}_1|N , D) = \chi_s(D) = -2 = \chi(F). \]

\begin{figure}
\labellist
\pinlabel $N$ at 217 110
\pinlabel \text{shift holonomy} at 135 -5
\pinlabel \text{shift holonomy} at 285 -5
\pinlabel \text{foliation transverse} at 0 15
\pinlabel {\text{to the} $[\frac{1}{4}, \frac{3}{4}]$ \text{factor}} at 0 0
\pinlabel {$S \times 1$} at 310 190
\pinlabel {$S \times \frac{3}{4}$} at 310 140
\pinlabel {$S \times \frac{1}{4}$} at 310 80
\pinlabel {$S \times 0$} at 319 37

\endlabellist

\centering
\includegraphics[width=3 in]{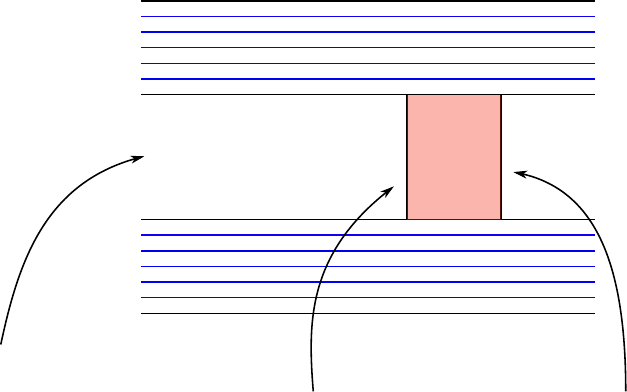}
\caption{A schematic picture of the constructed foliation for computing the Thurston norm. The horizontal foliations indicate the product foliation.}
\label{product-norm}
\end{figure}

Likewise, we can show 
\begin{eqnarray}
x([S]-[F])=-\chi(S)-\chi(F),
\label{thurston-norm-four}
\end{eqnarray}
by constructing a taut foliation $\mathcal{F}_2$ on $M$ whose Euler class $e(\mathcal{F}_2)$ assigns the numbers $\chi(S)$ and $-\chi(F)$ to the surfaces $S$ and $F$ respectively. We do the same steps except at the end we use a standard taut foliation of $N$ by a stack of monkey saddles whose relative Euler class assigns $-\chi_s(D)=2$ to the meridional disk of $N$. To sum up, we have proved Equalities (\ref{thurston-norm-one}), (\ref{thurston-norm-two}), (\ref{thurston-norm-three}), and (\ref{thurston-norm-four}). These four equalities show that the unit balls have the claimed shapes in Figure \ref{unitballs}. For example, Equalities (\ref{thurston-norm-one}), (\ref{thurston-norm-two}), and (\ref{thurston-norm-three}) determine the shape of the Thurston norm ball in the first and the third quadrants.

\end{proof}

\begin{lem}[Atoroidal and hyperbolic]
	Let $M$ be one of the manifolds constructed in Section \ref{construction-of-manifolds}. Then $M$ is irreducible, atoroidal and hyperbolic. 
	\label{atoroidal}
\end{lem}

\begin{proof}
	The proof is divided into two steps. In the first step, we show that any incompressible torus $T$ in $M$ can be isotoped to be disjoint from $N$. The second step is to prove that the manifold $M - N^\circ$ is atoroidal. Steps 1 and 2 imply that $M$ has no incompressible torus. On the other hand, the compact orientable 3-manifold $M \neq S^2 \times S^1$ admits taut foliations and hence is irreducible \cite{novikov1965topology,rosenberg1968foliations}. At this point, Thurston's hyperbolization theorem for Haken manifolds implies that $M$ is hyperbolic, since $M$ is closed, atoroidal and Haken.\\
	
	\textbf{Step 1}: Isotope $T$ such that $T \cap S $ is a collection of simple closed curves that are essential in both $T$ and $S$. This can be done since both $S$ and $T$ are incompressible and $M$ is irreducible. In particular, these curves are parallel in the torus $T$ and cut $T$ into a collection of annuli $T_i$ in $M \setminus \setminus S$. See Notation \ref{notation}. We may assume that $|T \cap S|$ is minimal. We show that $T_i$ can be isotoped to be disjoint from $N$. 
	
	Isotope $T_i$ such that the part of $T_i$ lying in $S \times [\frac{3}{4},1]$ and $S \times [0,\frac{1}{4}]$ becomes a standard vertical annulus. This can be done since $S \times [\frac{3}{4},1]$ and $S \times [0,\frac{1}{4}]$ admit product foliations and $T_i$ is essential, for example by Roussarie--Thurston procedure. Let $T_i'$ be the rest of $T_i$. Isotope $T_i'$ such that $T_i' \cap \partial N$ is a collection of simple closed curves in $\partial N$. These curves can be assumed to be essential in $\partial N$ and hence parallel to each other; otherwise one can reduce $|T \cap \partial N|$ after an isotopy of the essential surface $T_i'$. Let $A_1$ and $A_2$ be the annuli $\gamma_+ \times [\frac{1}{4}, \frac{3}{4}]$ and $\gamma_- \times [\frac{1}{4}, \frac{3}{4}]$ respectively. We may assume that the curves $T_i' \cap \partial N$ are not parallel to $\partial A_1$ and $\partial A_2$. This is because otherwise, each component of $T_i' \cap N$ is a $\partial$-parallel annulus in $N$ by Lemma \ref{incompressible-surface-in-torus}, and so it can be isotoped out of $N$. Hence, the intersection of $T_i'$ with $A_1$ and $A_2$ is a collection of properly embedded arcs in $A_1$ and $A_2$. Two types of arcs can occur: Type $1$ arcs start and end on the same component of $\partial A_1$ (or $\partial A_2$), whereas Type $2$ arcs run between different components of $\partial A_1$ (or $\partial A_2$). See Figure \ref{arcs-annulus}, left.\\

\begin{figure}
	\centering
	\begin{minipage}{4.2 cm}
	\labellist
	\small\hair 2pt
	\pinlabel {$\text{Type }1$}  at 325 150
	\pinlabel {$\text{Type }2$} at 175 100
	\pinlabel {$\gamma_+ \times \frac{3}{4}$} at 220 390
	\pinlabel {$\gamma_+ \times \frac{1}{4}$} at 220 -30
	
	\endlabellist
	
	\centering
	\includegraphics[width = 1.2 in]{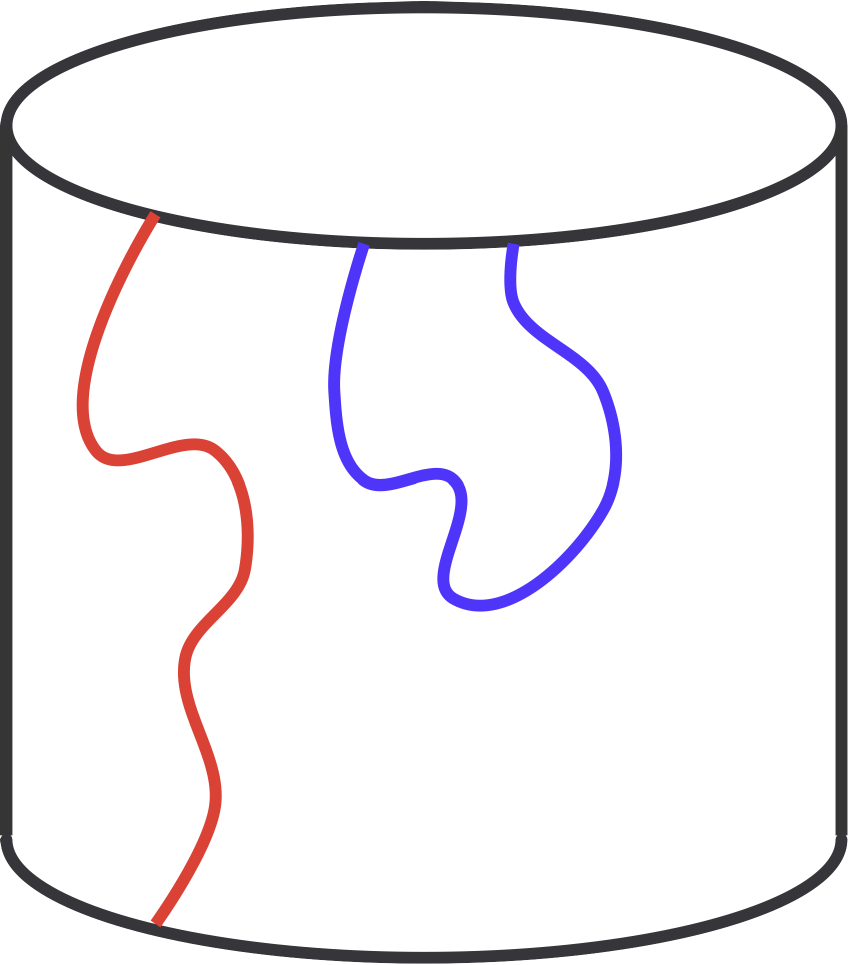}
	\end{minipage}
	\begin{minipage}{4.2 cm}
		\centering
		\includegraphics[width=1.2 in]{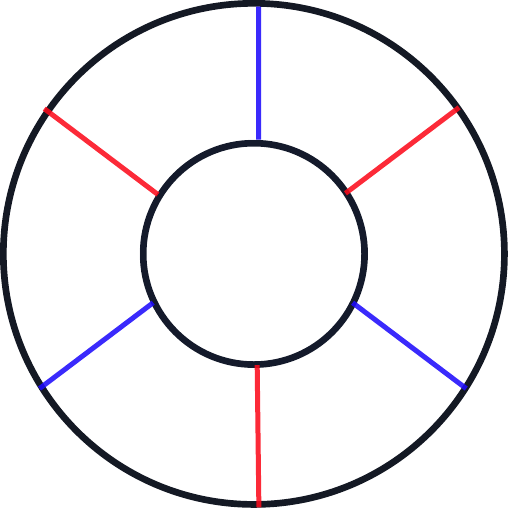}
	\end{minipage}
	\caption{Left: Two types of arcs in the annulus $A_1$. Right: Annulus $T_i'$ is cut into squares by Type 2 arcs.}
	\label{arcs-annulus}
\end{figure}
	
First, we remove Type 1 arcs. Choose an innermost Type $1$ arc on $A_1$ (or $A_2$). Assume that the endpoints of this arc lie on $S \times \{ \frac{3}{4}\}$. This arc together with a portion of the boundary of $A_1$ bounds a disk. Use this disk to push the arc, together with part of the corresponding annulus $T_i'$, out of $A_1$ and reduce the number of Type $1$ arcs. At this point, $T_i - T'_i$ might not be a vertical annulus any more, and we isotope $T_i$ further to make $T_i - T'_i$ a vertical annulus again. Denote by $\theta_1$ and $\theta_2$ the simple closed curves in $T'_i \cap (S \times \{ \frac{3}{4} \})$ containing the endpoints of the Type 1 arc; it is possible that $\theta_1=\theta_2$. From the perspective of $T'_i$, a neighborhood of the Type 1 arc is a band connecting $\theta_1$ to $\theta_2$, and the operation of pushing the Type 1 arc corresponds to surgering $\theta_1 \cup \theta_2$ along that band. This operation converts $\theta_1 \cup \theta_2$ into one or two simple closed curves one of which bounds a disk in $T_i$. This disk can be used to isotope $T_i$ and remove the non essential curve in $T_i$. In fact this argument shows that there must be two components after surgering along the band, otherwise $|T \cap S|$ would not have been minimal. At this point, $T_i - T'_i$ is a vertical annulus again. After doing this finitely many times, we are left with a number of Type $2$ arcs. 

From the perspective of $T_i'$, Type $2$ arcs start from one boundary component of $T_i'$ and end on the other boundary component of $T_i'$. See Figure \ref{arcs-annulus}, right. So $T_i' \cap N$ is a collection of squares. Each square in $T_i' \cap N$ has two sides on $\partial_{\pitchfork}N$ and two sides on $\partial_{\tau}N$ and hence is $\partial$-parallel in $N$. To see this, note that the meridional disk of $N$ has six sides on $\partial_{\pitchfork}N$ and six sides on $\partial_{\tau}N$. Starting from an innermost square, we can isotope them out of $N$. We have shown that $T_i$ can be isotoped to be disjoint from $N$. \\
	
	\textbf{Step 2}: We show that any incompressible torus in $M - N^\circ$ is $\partial$-parallel. Recall that $\gamma \subset S$ is an essential simple closed curve, and $A \subset S$ is a regular neighborhood of $\gamma$. Let $P: = M-N^\circ$. Then  $P= M_f -\text{int}\big(A \times [\frac{1}{4}, \frac{3}{4}]\big)$, where $M_f$ is the fibered hyperbolic 3-manifold with fiber $S$ and monodromy $f$ as described in \textbf{Step II} of Section \ref{construction-of-manifolds}. We may assume that $|T \cap S|$ is minimal.
	
	Define the map $\delta \colon \pi_1(P) \longrightarrow \mathbb{Z}$ as the algebraic intersection number with $S$, and denote by $h \colon \tilde{P} \longrightarrow P$ the infinite cyclic cover corresponding to the kernel of $\delta$. Then the lifts of $S$ in $\tilde{P}$ are indexed by elements of $\mathbb{Z}$, and we denote them by $S_j$ where $j \in \mathbb{Z}$. We can partition the annuli $T_i$ into $4$ groups corresponding to the way their lifts in $\tilde{P}$ look like, and we denote them by $X$, $Y$, $Z$ and $W$ as follows. If a lift of $T_i$ goes from $S_j$ to $S_{j'}$, then the index of $T_i$ is defined as $j'-j$. Here $X$ and $Y$ have indices $1$ and $-1$ respectively. The index of both $Z$ and $W$ is $0$, and they can be distinguished by looking at the side of $S$ that they lie in. See Figure \ref{index}.\\
	
	$P \setminus \setminus S$ has a natural foliation induced from the product foliation on $M_f \setminus \setminus S = S \times [0,1]$. Therefore, if $T_i$ has index $1$ or $-1$, then by Roussarie--Thurston procedure $T_i$ can be isotoped to be a vertical annulus. Note that technically $P \setminus \setminus S$ is not a sutured manifold, as it has concave corners rather than convex corners. Nevertheless since the foliation on $P\setminus \setminus S$ is almost a product, we can use the proof of Roussarie--Thurston rather than the statement to put $T$ in the desired position. Similarly, we use the Roussarie--Thurston procedure for an annulus $T_i$ of index $0$. In this case $T_i$ can be isotoped such that it consists of a union of tangential and transverse annuli. The tangential annuli lie inside a leaf of the foliation on $P \setminus \setminus S$, while transversal annuli have a product foliation induced from the foliation on $P \setminus \setminus S$.\\
	
	Denote the boundary components of $T_i$ by $c_i$ and $d_i$. Fixing a direction for $T$, we assume that $T_i$ starts from $c_i$ and ends at $d_i$. Hence, if the index of $T_i$ is $1$ (respectively $-1$), then $d_i$ is isotopic to $f (c_i)$ (respectively $f^{-1} (c_i)$). Similarly if $T_i$ has index $0$, then $d_i$ is isotopic to $c_i$. Define the total index, $I$, of $T$ as sum of the indices of all $T_i$. It follows that $f^I(c_1)$ is isotopic to $c_1$. As $f$ is pseudo-Anosov and $c_1 \subset S$ is an essential closed curve, $I$ should be equal to $0$.
	
	For any $T_i$ of index $0$, let $D_i$ be a $\partial$-compressing disk for $T_i$ inside $S \times [0,1]$. Then one can use $D_i$ and isotope $T_i$ across $S$ to reduce $|T \cap S|$ unless $D_i$ intersects the core of the solid torus $A \times [\frac{1}{4}, \frac{3}{4}]$ essentially, that is algebraically a nonzero number of times.  Note that this core curve can be identified with $\gamma \times \{ \frac{1}{2} \}$. In this case, $c_i$ and $d_i$ are both isotopic to $\gamma$ in $S$ if $T_i$ is of type $Z$, and both are isotopic to $f(\gamma)$ if $T_i$ is of type $W$. Call this condition $(*)$.\\
	
	Now we show that $|T \cap S|=0$ and $T$ is $\partial$-parallel in $M - N^\circ$. Assume to the contrary. Note that there should be at least one annulus of index $0$, since otherwise the total index is nonzero. If we allow both annuli of index $0$ and of index nonzero, then the condition $(*)$ would be violated for the following reason. It is clear that looking at the sequence of annuli in $T$, right after an $X$ only $X$ or $Z$ can appear. Similarly right after a $Y$ there can come a $Y$ or $W$. Right after $Z$ comes a $Y$ or $W$, and right after $W$ only $X$ or $Z$ can happen. Now suppose letters $a \in \{X, Y \}$ and $b \in \{ Z, W \}$ appear in the sequence corresponding to $T$. It follows that there should be a sequence of consecutive letters that starts and ends with elements in $\{ Z , W \}$, while all the intermediate letters belong to $\{ X , Y \}$. It is easy to see that the only possibilities are the sequences $Z Y \cdots YW$ and $W X \cdots XZ$. But in both cases, the condition $(*)$ is violated by at least one of the two ending letters $W$ and $Z$. This is because the intermediate letters change the curve $c_i$ by powers of $f$, and the existence of such string of letters implies that $f^r(\gamma)$ is isotopic to $\gamma$ for some $r \in \mathbb{N}$, contradiction to $f$ being pseudo-Anosov. So all the annuli $T_i$ are of index $0$, and there are at least two annuli $T_i$. But again the condition $(*)$ is violated. The contradiction shows that the torus $T$ is disjoint from $S$. \\
	
	\begin{figure}
		\labellist
		\pinlabel $S_{j+1}$ at 20 110
		\pinlabel $S_j$ at 300 10
		\endlabellist
		
		\centering
		\includegraphics[width=2.5 in]{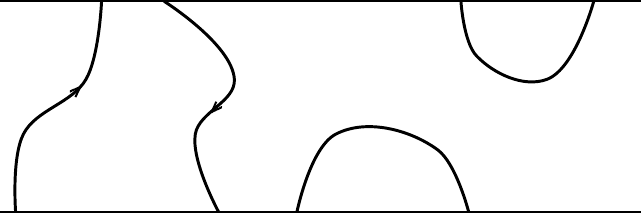}
		\caption{The annuli $T_i$ are divided into four types according to the way their lifts in the covering space $\tilde{P}$ look like. The pictures from left to right respectively correspond to the lift of an annulus of type $X$, $Y$, $Z$, and $W$.}
		\label{index}
	\end{figure}

	Use Roussarie--Thurston procedure to isotope $T$ inside the foliated manifold $P \setminus \setminus S$ so that $T$ consists of tangential and transverse annuli. The cross section of the picture is a simple closed curve winding around a point, where the point is a cross section of $\gamma \times \{ \frac{1}{2} \}$. Since the curve is simple and closed, the winding number, $w$, is $0$ or $\pm 1$. If $w = 0$, then $T$ is not incompressible. Therefore, $w = 1$ or $-1$, which correspond to a $\partial$-parallel torus in $P$.
\end{proof}

For the proof of the next lemma, we need the following less known theorem of Gabai. 

\begin{thm}[Gabai \cite{gabai1983foliations}]
	Let $M$ be a closed, orientable 3-manifold. Assume that $P$ and $Q$ are two possibly disconnected, norm-minimizing surfaces in $M$ that are homologous. There is a sequence of possibly disconnected, norm-minimizing surfaces $P=P_0 , P_1 , \cdots , P_n=Q$ with each term in the same homology class as $[P]=[Q]$ such that any two adjacent terms in the sequence can be isotoped to be disjoint in $M$.
	\label{disjointsurfaces}
\end{thm}

\begin{proof}
	Directly follows from the proof of Lemma 3.6. in \cite{gabai1983foliations}.
\end{proof}

\begin{lem}
	Let $M$ and $S$ be as defined in Section \ref{construction-of-manifolds}. Any norm-minimizing surface that is homologous to $S$ is in fact isotopic to $S$.
	\label{unique}
\end{lem}

\begin{proof} Suppose $S'$ is a norm-minimizing surface that is homologous to $S$. By Theorem \ref{disjointsurfaces}, we can assume that $S'$ is disjoint from $S$, and hence lies inside $M \setminus \setminus S$. Recall that $\gamma$ is a non-separating simple closed curve inside $S$, and $A$ is a tubular neighborhood of $\gamma$ is $S$ with $\partial A = \gamma_+ \cup \gamma_-$. The manifold $M \setminus \setminus S$ is the union of two parts: the \emph{product part} $(S- A^\circ) \times [0,1]$ and the \emph{twisted part} 
\begin{eqnarray}
 N' := (A \times [0,\frac{1}{4}]) \cup N \cup (A \times [\frac{3}{4},1]). 
 \label{N'}
\end{eqnarray}

\begin{figure}
	\labellist
	\pinlabel $N$ at 150 50 
	\pinlabel $N'$ at 412 50
	\endlabellist
	
	\centering
	\includegraphics[width=4 in]{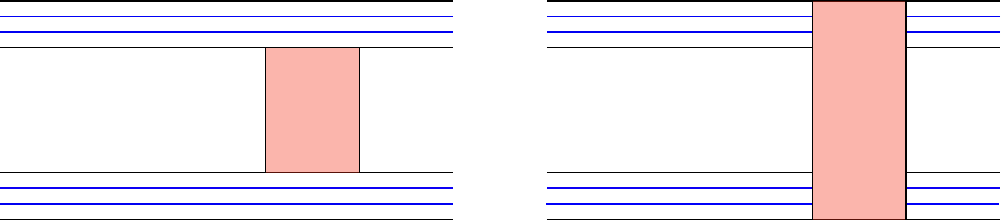}
	\caption{Submanifolds $N$ and $N'$ of $M \setminus \setminus S$. The top and bottom parts with the product foliation on the left indicate $S \times [\frac{3}{4}, 1]$ and $S \times [0, \frac{1}{4}]$.}
	\label{N-and-N'}
\end{figure}	
See Figure \ref{N-and-N'}. Then $N'$ is a sutured manifold with annuli sutures $A'_1 : = \gamma_+ \times [0,1]$ and $A'_2 := \gamma_- \times [0,1]$. Note that $N'$ and $N$ are homeomorphic as sutured manifolds, since $N'$ is obtained by thickening the tangential boundary of $N$ to the outside.  We isotope $S'$ to be in a standard form in each of the two parts of $M \setminus \setminus S$ and conclude that $S'$ is isotopic to $S$.

By assumption, $S'$ is disjoint from the top and bottom copies of $S$ in $M \setminus \setminus S$. Isotope $S'$ such that it intersects $A'_i$ in a union of simple closed curves. Inessential circles in $S' \cap A'_i \subset A'_i$ can be removed, starting from an innermost circle and isotoping the incompressible $S'$. At this point $S' \cap A'_i$ is a union of disjoint copies of essential circles, each of them is isotopic to the core curve of $A'_i$. Let 
\[ S'_1 := S' \cap \Big((S-A^\circ) \times [0,1] \Big) , \hspace{6mm}  S'_2 := S' \cap N'. \]
After further isotopy of $S'$, it can be assumed that no component of $S'_1$ is parallel to a subannulus in $A'_1$ or $A'_2$. We claim that $S'_1$ does not admit any $\partial$-compressing disk $D_1$ with $\partial D_1 = \alpha_1 \cup \beta_1$, where $D_1 \cap S'_1 = \alpha_1$ is an essential arc, $\beta_1 \subset A'_1 \cup A'_2$ is an arc, and $\partial \alpha_1 = \partial \beta_1 = \alpha_1 \cap \beta_1$. Assume to the contrary that such a $\partial$-compressing disk $D_1$ exists. Then the two points in $\partial \alpha_1$ lie on distinct components $b_1$ and $b_2$ of $\partial S'_1 \subset A'_1 \cup A'_2$. Let $\hat{A}$ be the annulus cobounding $b_1$ and $b_2$ in $A'_1 \cup A'_2$.  Let $N_1$ be a neighborhood of $b_1 \cup b_2 \cup \alpha_1$ in $S'_1$, which is a 3-punctured sphere. The circle $\partial N_1 - \partial S'_1$ bounds a disk in the complement of $S'_1$, lying near $D_1 \cup \hat{A}$. Since $S'_1$ is incompressible this boundary circle also bounds a disk in $S'_1$. Thus a component of $S'_1$ is an annulus, which we refer to by $S'_1$ by abuse of notation. Surgering the torus $S'_1 \cup \hat{A}$ along $D_1$ yields a sphere, which bounds a ball since $(S-A^\circ) \times [0,1]$ is irreducible. Hence some component of $S'_1$ is a $\partial$-parallel annulus, contradicting our previous assumption. The contradiction shows that such $\partial$-compressing disk does not exist.

Now Lemma \ref{incompressible-product} shows that $S'_1$ is a union of parallel copies of $S- A^\circ$. But there cannot be more than one copy, since otherwise 
\[ |\chi(S')| \geq 2 |\chi(S-A^\circ)| = 2 \times (2g-2) > 2g-2 = |\chi(S)|,  \]
contradicting the assumption that $S'$ is norm-minimizing and $[S']=[S]$. Hence $S'_1$ is a single copy of $S-A^\circ$, and $S'$ intersects each of $A'_i$ in exactly one essential simple closed curve. But the only incompressible surfaces in $(N', A_1' \cup A_2')$ with this oriented boundary are the two components of the tangential boundary of $N'$, up to isotopy (Lemma \ref{incompressible-surface-in-torus}). Therefore, $S'$ is isotopic to either the top or the bottom copy of $S$ in $M \setminus \setminus S$. This completes the proof. 
\end{proof}

\begin{lem}[General position for annuli] 
Let $M$ be one of the manifolds constructed in Section \ref{construction-of-manifolds}, and $\mathcal{F}$ be a taut foliation on $M \setminus \setminus S$ having $\partial (M \setminus \setminus S)$ as leaves. Recall that 
\[ M \setminus \setminus S \cong \Big( (S-A^\circ)\times[0,1] \hspace{2mm} \Big) \cup \hspace{2mm} N', \]
where $N'$ is the solid torus defined as in Equation (\ref{N'}), and $\partial (S- A^\circ) \times [0,1] = A'_1 \cup A'_2$. Then $A'_1$ and $A'_2$ can be isotoped relative to their boundaries such that the induced foliations on them are suspension foliations. 
\label{annuli-general-position}
\end{lem}

\begin{proof}
	Let $K := S-A^\circ$. First we show that $A'_1 \cup A'_2$ can be isotoped relative to boundary such that $\mathcal{F}$ is transverse to $ A'_1 \cup A'_2$ and there is no annulus leaf in $(K \times [0,1], A'_1 \cup A'_2)$. After such an isotopy, the induced foliation on $K \times [0,1]$ contains no Reeb component or half Reeb component. 
	
	By Roussarie--Thurston general position, the incompressible $A'_1 \cup A'_2$ can be isotoped to be transverse to the foliation, although there might be two-dimensional Reeb components on $A'_i$. Consider the induced foliation on $K \times [0,1]$, and by abuse of notation call it $\mathcal{F}$ again. Leaves of $\mathcal{F}$ are incompressible by Novikov's theorem. 
	
	We show that any annulus leaf $L$ in $K \times [0,1]$ is parallel to an essential subannulus in $A'_1$ or $A'_2$. If the annulus $L$ admits no $\partial$-compressing disk $D$ with $\partial D = \alpha \cup \beta$, $D \cap L = \alpha$ an essential arc, $\beta \subset \partial K \times [0,1]$ an arc, and $\alpha \cap \beta = \partial \alpha= \partial \beta$, then $L$ would be horizontal by Lemma \ref{incompressible-product}, which is not possible. The contradiction shows that such a $\partial$-compressing disk $D$ exists. Without loss of generality assume that $\partial L \subset A'_1$ and let $\hat{A}$ be the portion of $A'_1$ cobounding the components of $\partial L$. Let $P$ be the torus $L \cup \hat{A}$. Then $P$ surgered along $D$ is a sphere that bounds a solid ball since $K \times [0,1]$ is irreducible. This ball can be used to isotope $L$ into $\hat{A}$. Hence every annulus leaf is parallel to a subannulus in $A'_i$. 
	
	By Haefliger, there is an outermost annulus leaf parallel to an essential subannulus of $A'_i$. Therefore, we may isotope $A'_i$ to push all annuli leaves out of $K \times [0,1]$. After this isotopy, the foliation is still transverse to $A'_i$, and there is no annulus leaf in $(K \times [0,1], A'_1 \cup A'_2)$.\\
	
	Now the induced foliation on $K \times [0,1]$ is an \emph{essential lamination} on an $I$-bundle over the compact surface $K$ of negative Euler characteristic, and having $K \times \{ 0 , 1 \}$ as leaves. By Brittenham \cite{brittenham1999essential, brittenham1997essential}, this essential lamination can be isotoped to be transverse to $[0,1]$ factor. In particular, after the isotopy there are no Reeb components on $A'_i$.
\end{proof}

\section{Taut foliations on sutured solid tori}
\label{section:sutured-solid-tori}

\begin{definition}
	A taut foliation $\mathcal{F}$ of a sutured solid torus $(T, \gamma)$ is \emph{standard} if either 
	\begin{enumerate}
		\item $\mathcal{F}$ is obtained from a stack of generalized saddles by possibly $I$-bundle replacement along some of the components of its tangential boundary (see Example \ref{stackexample}), or
		\item $\mathcal{F}$ is a foliation of $\text{annulus} \times I$, transverse to the $I$ factor, where $\gamma = \partial (\text{annulus}) \times  I$.
	\end{enumerate} 
	\label{standard}
\end{definition}

In the special case above, the $I$-bundle replacement is the operation of replacing a compact leaf $L$ with an $I$-bundle over $L$ which is foliated transverse to $I$ factor ($I$ an interval). 

\begin{definition}
	Let $(M, \gamma)$ be a sutured manifold, and $L$ be a properly embedded surface in $M$. We say $L$ is \emph{$\partial_\tau M$-parallel} if there is a union $R$ of components of $R(\gamma)$ such that $R$ and $L$ together bound a submanifold homeomorphic to $R \times [0,1]$, where $R \times 0$ (respectively $R \times 1$) is identified with $R$ (respectively $L$), and $\partial R \times [0,1] \subset \gamma$. 
\end{definition}

In this section, we prove Proposition \ref{maximal euler class} and Corollary \ref{nofoliation}. The following definition is essentially due to Novikov \cite[Page 6]{novikov1965topology}.

\begin{definition}
Let $\mathcal{F}$ be a transversely oriented foliation on a compact 3-manifold $M$ and $L$ be a leaf of $\mathcal{F}$. Fix a point $p \in L$, and let $\mathcal{T}(L,p)$ be the set of positive closed transversals for $L$ that start and end at $p$. Positive means that its orientation agrees with the transverse orientation of $\mathcal{F}$. The set of based homotopy classes of elements of $\mathcal{T}(L,p)$ forms a semigroup under concatenation. We call it the \emph{based transversal semigroup of $L$} and denote it by $\mathcal{T}_{h}(L,p)$. 
\label{transversal-semigroup}
\end{definition}

 \begin{remark}
 	\begin{enumerate}[i)]
 		\item For different choices $p$ and $q$ of base points, the semigroups $\mathcal{T}_h(L,p)$ and $\mathcal{T}_h(L,q)$ are isomorphic; the isomorphism depends only on the homotopy class on the leaf $L$ of a path joining $p$ and $q$ on $L$. To see this, pick an oriented arc $\delta \subset L$ from $q$ to $p$. Let $\gamma$ be a positive transversal based at $p$ and denote by $\hat{\gamma}$ the concatenation of $\delta$, $\gamma$ and $-\delta$. A perturbation of $\hat{\gamma}$ is a positive transversal based at $q$. See Novikov \cite[Lemma 2.2]{novikov1965topology}.
 		\item In general the set $\mathcal{T}(L,p)$ might be empty. However, when $\mathcal{F}$ is taut and $M$ is closed, $\mathcal{T}(L,p)$ is non-empty for each leaf $L$ and $p \in L$ by definition. 
 		\item If $L$ is noncompact and $\mathcal{F}$ is not necessarily taut, then $\mathcal{T}(L,p)$ is non-empty. This is because the compact manifold can be covered by finitely many foliation charts, and hence some chart intersects the noncompact leaf $L$ infinitely many times. Take a small transverse arc in the chart with endpoints on $L$ and connect its endpoints by a path in $L$ to obtain a closed curve. This closed curve can be perturbed to a transversal which is possibly not based at $p$. Now the base point can be moved to $p$ while keeping the curve a positive transversal. 
 		\item By Novikov's theorem, for Reebless $\mathcal{F}$ every element in $\mathcal{T}(L,p)$ is homotopically nontrivial.
 		\item In \cite[Page 6]{novikov1965topology}, Novikov defined a semigroup very similar to $\mathcal{T}_h(L,p)$; he considered regular homotopy classes of transversals where the initial point is always at $p$ and during the homotopy the transversal is never tangent to leaves. The semigroup defined above is a quotient of Novikov's semigroup.
 	\end{enumerate}
 \end{remark}

\begin{figure}
\centering
\includegraphics[width = 3.5 in]{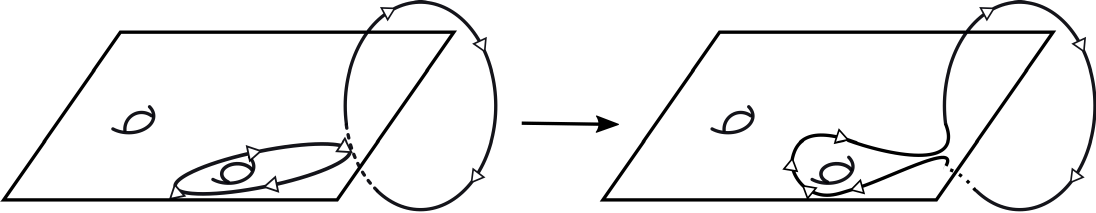}
\caption{The concatenation of a based positive transversal and a based loop in $L$ can be perturbed to a new based positive transversal.}
\label{transloop}
\end{figure}

\begin{obs} 
If $\mathcal{T}_h(L,p)$ is non-empty then it is closed under multiplication by elements of $i_*(\pi_1(L,p))$, where $i:L \longrightarrow M$ is the inclusion map.
\label{product}
\end{obs}
\begin{proof}
See Figure \ref{transloop}.
\end{proof}

%


\begin{cor}
Let $(T,\gamma)$ be a sutured solid torus, and $\mathcal{F}$ be a taut foliation on $T$. Any noncompact leaf of $\mathcal{F}$ is simply connected. 
\label{simplyconnected}
\end{cor}

\begin{proof} Assume to the contrary. Let $L$ be a noncompact leaf such that $\pi_1(L, p) \neq \{ 0 \}$ where $p \in L$. By Novikov's theorem, $L$ is $\pi_1$-injective. Hence if $l$ is a generator of $\pi_1(T,p)$, then $i_*(\pi_1(L,p)) = \langle l^k \rangle$ for some natural number $k$. Since $L$ is noncompact, $\mathcal{T}(L,p)$ is non-empty. Let $\gamma \in \mathcal{T}(L,p)$ be a transversal, and assume that $\gamma$ represents $l^m \in \pi_1(T, p)$ for some $m \in \mathbb{Z} \setminus \{ 0 \}$. The iterated concatenation $\gamma^k$ is also a based transversal that represents the element $l^{km} \in \mathcal{T}_h(L,p)$. By Observation \ref{product}, we have $0 = (l^{-k})^{m} \cdot l^{km} \in \mathcal{T}_h(L,p)$, since $l^{-k} \in i_*(\pi_1(L,p))$. This contradicts Novikov's theorem. Therefore $L$ is simply connected. 
\end{proof}

\begin{prop}
	Let $(T,\gamma)$ be a sutured solid torus, and $\mathcal{F}$ be a taut foliation of $T$ such that every compact leaf of $\mathcal{F}$ is $\partial_\tau T$-parallel. Then $\mathcal{F}$ is standard. In particular, $\mathcal{F}$ has maximal relative Euler class, meaning that its pairing with the meridional disk $D$ is equal to $\chi_s(D)$ up to sign. Here $\chi_s(S)$ denotes the sutured Euler characteristic of $D$, where $D$ is assumed to intersect $\gamma$ minimally. 
	\label{maximal euler class}
\end{prop}

\begin{proof}
	
First, it can be assumed that every compact leaf of $\mathcal{F}$ is a component of $R(\gamma)$. By hypothesis, every compact leaf of $\mathcal{F}$ is a $\partial _{\tau} T$-parallel annulus. So by Haefliger's theorem, for every component $S$ of $R(\gamma)$ there is a compact leaf $L$ such that the foliation between $L$ and $S$ is a suspension foliation of $\text{annulus} \times \text{interval}$ transverse to the $\text{interval}$ factor, and $L$ is farthest from $S$; we include the possibility of $L = S$. Remove the foliation between $L$ and $S$, and repeat this for all components of $R(\gamma)$. By abuse of notation, call the new foliation $\mathcal{F}$. Assume that $\mathcal{F}$ is non-empty, otherwise the conclusion holds.

By Corollary \ref{simplyconnected}, \textit{every noncompact leaf of $\mathcal{F}$ is simply connected}. So, \textit{if} $\mu \in \text{Homeo}^+{(I)}$ \textit{is the holonomy of an annulus suture, then $\mu$ has no fixed point except the interval endpoints}; otherwise a fixed point corresponds to a nontrivial loop inside a noncompact leaf. We want to prove that the foliation is indeed the standard foliation by a stack of generalized saddles with holonomy $\mu$. A small foliated neighborhood of $\partial_{\tau}T$ is determined completely by the holonomy $\mu$; see e.g. Camacho and Neto \cite[Page 67]{camacho2013geometric}. See Example \ref{stackexample}. Since $\mu$ is a shift, near $\partial_{\tau}T$ the leaves spiral around $\partial_{\tau}T$ in the standard way, i.e. like the picture for a stack of generalized saddles with the shift holonomy. Hence if we shave a small neighborhood of the tangential boundary, we obtain a foliation, $\mathcal{G}$, on a solid torus $T'$ that is transverse to the boundary and whose picture looks like Figure \ref{shave}. Note that $\mathcal{G}$ is a subfoliation of $\mathcal{F}$ and therefore has no Reeb component. Pick a curve $s \subset \partial T'$ that is isotopic to the sutures and is transverse to the foliation $\partial \mathcal{G}$. Since the foliation $\partial \mathcal G$ is transverse to $s$, it has no Reeb component and therefore is a suspension foliation. We prove that the foliation $\mathcal{G}$ is a product foliation by disks. 

Let $m \subset \partial T'$ be a meridian. Since the foliation $\partial \mathcal{G}$ is Reebless, after an isotopy either $m$ is transverse to $\partial \mathcal{G}$ or $m$ is a leaf of $\partial \mathcal{G}$. By Novikov, $m$ cannot be transverse to the Reebless $\mathcal{G}$, and so $m$ is a leaf of $\partial \mathcal{G}$. Let $Q$ be the leaf of $\mathcal{G}$ with $m \subset \partial Q$. By Novikov, $Q$ is $\pi_1$-injective. Since $m$ bounds a disk in $T'$ and $Q$ is $\pi_1$-injective, $Q$ itself should be a disk. Now by the Reeb stability theorem, $\mathcal{G}$ is the product foliation by disks. 

Since $\mathcal{F}$ was obtained by adding standard pieces to the boundary of $\mathcal{G}$, $\mathcal{F}$ is a standard taut foliation by a stack of generalized saddles, and the relative Euler class of $\mathcal{F}$ is maximal.

\begin{figure}
\centering
\includegraphics[width=3 in]{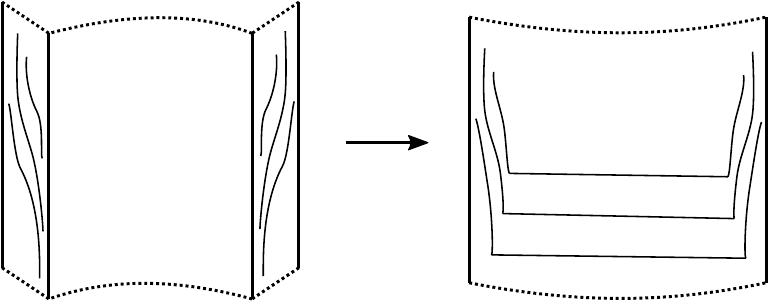}
\caption{Shaving a neighborhood of $\partial_{\tau}T$}
\label{shave}
\end{figure}




\end{proof}

\begin{cor}
	Let $(T, \gamma)$ be a sutured solid torus with either two or four sutures. Every taut foliation of $T$ is standard. In particular, the relative Euler class of the foliation is maximal, meaning that it assigns $\chi_s(D)$ to the meridional disk $D$ of $T$ up to sign.
	\label{nofoliation} 
\end{cor}

\begin{proof}
By Proposition \ref{maximal euler class}, it is enough to show that every compact leaf $A$ is a $\partial_{\tau}T$-parallel annulus. By Novikov, any compact leaf $A$ is $\pi_1$-injective, and hence is a $\partial$-parallel annulus or a disk (Lemma \ref{incompressible-surface-in-torus}). No leaf is a disk either, otherwise by Reeb Stability theorem the foliation should be the product $D^2 \times S^1$. Note that $(A, \partial A) \subset (T, \gamma)$. Moreover, each component of $\partial A$ is parallel to the core curve of the corresponding suture that it lies on, and their orientations agree. Since the number of sutures is either $2$ or $4$, (at least) one of the two components of $\partial T - \partial A$ contains exactly one tangential annulus. Therefore, $A$ is $\partial_{\tau}T$-parallel.
\end{proof}

\begin{remark}
Example \ref{stackexample} shows that Corollary \ref{nofoliation} is not true when there are more than four sutures.
\end{remark}

\begin{remark}
	Given a sutured solid torus, one can classify taut foliations on it using the following steps: 
	\begin{enumerate}[i)]
		\item Locate the compact annuli leaves.
		\item Excise foliated packets containing all compact leaves, and produce a union of sutured solid tori such that every compact leaf is a component of the tangential boundary. 
		\item By Proposition \ref{maximal euler class}, the foliation on each new sutured solid torus is a standard one by a stack of generalized saddles. Now glue the sutured solid tori pieces back together. 
	\end{enumerate} 
\end{remark}

\section{Proof of the main theorem}
	
	We prove that the manifold $M$ constructed in Section \ref{construction-of-manifolds} violates Thurston's conjecture. Consider the cohomology class $a = (0 , 2-2g) \in 2H^2(M ; \mathbb{Z})$ in Figure \ref{unitballs}, right. This class has dual Thurston norm equal to one (Lemma \ref{unitballs}). We show that $a$ cannot be realized as the Euler class of any taut foliation on $M$. Assume to the contrary, that there is such a taut foliation $\mathcal{F}$.\\  
	
	\textbf{Step 1}: We show that it can be assumed that $S$ is a leaf. Note that $S$ is an algebraically fully marked surface. Moreover, $M$ is hyperbolic (Lemma \ref{atoroidal}) and the surface $S$ is the unique norm-minimizing surface in its homology class, up to isotopy (Lemma \ref{unique}). By the fully marked surface theorem, there exists a new taut foliation that has $S$ as a leaf and whose oriented tangent plane field is homotopic to that of $\mathcal{F}$. In particular, the new foliation has the same Euler class as $\mathcal{F}$. By abuse of notation, we call the new foliation $\mathcal{F}$ again.\\
	
	\textbf{Step 2}: Recall that $A$ was an annulus neighborhood of $\gamma$. Note that if $M_1 = M \setminus \setminus S$, then there exists an embedded $(S-A^\circ) \times I$ in $M_1$ where $(S-A^\circ) \times \{ 0 \} $ and $(S-A^\circ) \times \{ 1 \} $ are in $\partial M_1$, and $\partial (S-A^\circ) \times I$ are vertical annuli. See Figure \ref{N-and-N'}. Apply Lemma \ref{annuli-general-position} and isotope $A'_1 \cup A'_2 = \partial (S-A^\circ)  \times I$ relative to their boundaries such that the induced foliations on $A'_1$ and $A'_2$ are suspension foliations. Note $N' = M_1 - \text{int}\big((S-A) \times I\big)$ is homeomorphic to $N$ as sutured manifolds. See Figures \ref{solidtorus} and \ref{N-and-N'}. The induced foliation on the sutured solid torus $N'$ is taut. This is because $\mathcal{F}_1$ is taut and hence every leaf of $N$ has either a closed transversal, or a transverse arc going between the components of $\partial (M \setminus \setminus S)$. Now given that the induced foliations on $A'_1 \cup A'_2$ are suspension foliations, it is easy to see that every leaf of $N'$ has either a closed transversal or a transverse arc connecting $A \times 0$ to $A \times 1$, i.e. $\mathcal{F}'$ is taut.\\
	
	\textbf{Step 3}: Let $D$ and $D'$ be the meridional disks of $N$ and $N'$ respectively, and $F$ be the surface of genus two constructed in the proof of Lemma \ref{Thurston-norm}. Recall that $F$ is obtained from $D$ by attaching the three bands $\delta_i \times [\frac{1}{4},\frac{3}{4}]$ and the two annuli $\alpha \times [\frac{3}{4},1]$ and $\beta \times [0, \frac{1}{4}]$ and then identifying $\alpha \times 1$ with $\beta \times 0$. Paraphrasing this, $F$ is obtained from $D'$ by attaching the three bands $\delta_i \times [0,1]$ and then identifying $\alpha \times 1$ with $\beta \times 0$. 
	
	Let $\mathcal{F}_1$ and $\mathcal{F}'$ be the induced foliations on $M_1$ and $N'$ respectively. Both $\mathcal{F}_1$ and $\mathcal{F}'$ are taut. Let $M_2 = (S-A^\circ) \times [0,1]$, and denote the induced foliation on $M_2$ by $\mathcal{F}_2$. Isotope the three bands $\delta_i \times [0,1]$ in $M_2$ relative to their boundaries such that the induced foliations on them are products. More precisely, by Roussarie--Thurston general position for the sutured manifold $M_2$ and taut foliation $\mathcal{F}_2$, the bands can be isotoped relative to their boundaries such that the induced foliations on them have only saddle singularities. By the Poincar\'{e}--Hopf formula, the number of such saddle tangencies is equal to zero. On the other hand since a band is simply connected, there is no room for a Reeb component either. Hence, the induced foliations on the bands are product foliations.\\

	\textbf{Step 4}: Let $F_1 \subset M_1$ be the surface obtained by cutting $F$ along $S$; in other words $F_1$ is the union of $D'$ and the three bands $\delta_i \times [0,1]$. We have
	\[ 0 = \langle e(\mathcal{F}), [F] \rangle = \text{Ind}(\mathcal{F}, F)= \text{Ind}(\mathcal{F}_1, F_1) = \text{Ind}(\mathcal{F}', D').  \]
	Here the first equality is by the assumption on the initial Euler class. The second and third equalities follow from the index sum formula, since the induced foliations on the bands are products and contain no singularity. Thus we have a taut foliation on the sutured manifold $N' \cong N$ with relative Euler class zero. However, by Corollary \ref{nofoliation} such a taut foliation does not exist. This gives a contradiction, and completes the proof of the main theorem.  


\section{Euler classes of general foliations}

For completeness and comparison, we bring the following result of Wood \cite[Page 351--352]{wood1969foliations}.

\label{generalfoliations}
\begin{thm}[Wood]
Let $M$ be a closed orientable 3-manifold and $a \in 2 H^2(M ; \mathbb{Z})$, i.e. the cohomology class $a$ satisfies the parity condition. There is a transversely oriented (not necessarily taut) foliation $\mathcal{F}$ of $M$ such that $e(\mathcal{F}) = a$.
\label{Wood}
\end{thm}

\begin{proof}
By Wood, every transversely oriented plane field on $M$ is homotopic to an integrable plane field, i.e. coming from a foliation \cite{wood1969foliations}. Homotopic plane fields have the same Euler class. Therefore, it suffices to find a plane field over $M$ with Euler class equal to $a$. 

Fix a trivialization of the tangent bundle $\mathrm{T}M$ of $M$, and identify $\mathrm{T}M \cong M \times \mathbb{R}^3$. A transversely oriented plane fields $\sigma$ over $M$ defines a map $f \colon M \longrightarrow S^2$, where $f(m)$ for $m \in M$ is the oriented unit normal vector to $\sigma$. Note that $\sigma = f^*(\mathrm{T}S^2)$, where $\mathrm{T}S^2$ is the tangent bundle to the 2-sphere. Therefore, if $s$ is the positive generator for $H^2(S^2 ; \mathbb{Z})$, then 
\[ e(\sigma) = f^*(e(\mathrm{T}S^2)) = f^*(2s)= 2 f^*(s).  \]
So it is enough to show that every element in $H^2(M ; \mathbb{Z})$ can be obtained by pulling back the fundamental cohomology class of $S^2$ under some map $f \colon M \longrightarrow S^2$. A class $a \in H^2(M; \mathbb{Z})$ corresponds to the homotopy class of a map $f \colon M \longrightarrow K(\mathbb{Z}, 2)= \mathbb{C}P^\infty$ which factors through the 2-skeleton $S^2$ up to homotopy. This completes the proof.
\end{proof} 

\section{Further Questions}
\label{section:questions}
\textbf{Taut foliations with trivial Euler class, and Anosov flows}: 
Although Thurston's conjecture was stated for the unit sphere of dual norm, a priori we only know that the Euler class of a taut foliation is inside the unit ball. 

\begin{question}
Which points inside the dual unit ball can be realized as the Euler class of some taut foliation on $M$? 
\end{question}

A point strictly inside the dual unit ball cannot correspond to a taut foliation having a compact leaf of negative Euler characteristic. This makes it difficult to construct taut foliations with Euler class strictly inside the dual unit ball. An interesting case is that of taut foliations with trivial Euler class. \emph{Anosov flows} provide one way of constructing taut foliations of Euler class zero. The \emph{weak stable (unstable) foliation} of an Anosov flow is a taut foliation since it has no compact leaves, and it has trivial Euler class since the flow direction is a section. To the best of my knowledge, it is not known whether or not every closed hyperbolic 3-manifold with positive first Betti number admits an Anosov flow. 
\begin{question}
Which 3-manifolds with positive first Betti number admit a taut foliation with trivial Euler class? 
\end{question}
As a concrete example, we do not know if the Whitehead link complement has a taut foliation, transverse to the boundary, and with relative Euler class zero. Such a taut foliation, if it exists, should have boundary Reeb components because of the \emph{relative parity condition}. More precisely, if there is no Reeb component on $\partial M$, then $\chi_s(S)$ and $\langle e(\mathcal{F}) , [S] \rangle$ have the same parity for each properly embedded oriented surface $S$. Now if $\mathcal{F}$ has no boundary Reeb component, and $S$ is a twice-punctured disk bounding one of the link components, then 
\[  \langle e(\mathcal{F}),[S] \rangle \equiv \chi_s(S) \equiv 1 \hspace{3mm} (\text{mod }2), \]
implying that the relative Euler class $e(\mathcal{F})$ is nonzero. \\


\textbf{Pseudo-Anosov and quasigeodesic flows}:
This subsection contains no new questions; instead it includes a brief and non exhaustive summary of what is known about realization of Euler classes of pseudo-Anosov or quasigeodesic flows. Cannon and Thurston \cite{cannon2007group} proved that the suspension flow of a pseudo-Anosov automorphism of a surface can be chosen to be quasigeodesic and pseudo-Anosov. Mosher, following Gabai, generalized \cite{cannon2007group} and proved that for any finite depth taut foliation $\mathcal{F}$ on a closed, hyperbolic 3-manifold there exists an \emph{almost transverse} pseudo-Anosov flow \cite{mosher1examples,mosher1996laminations}. Almost transverse means that it will be transverse to $\mathcal{F}$ after an appropriate blow up of a finite collection of closed orbits \cite{mosher1996laminations}. Fenley and Mosher proved that these flows are quasigeodesic as well \cite{fenley2001quasigeodesic}. Since the foliations constructed by Gabai in Theorem \ref{gabai} are of finite depth, it follows that the vertices of the dual ball are realized as Euler classes of pseudo-Anosov (quasigeodesic) flows on closed hyperbolic 3-manifolds.\\

\textbf{Euler classes of representations into $\Text{Homeo}^+(S^1)$}:
The following question is inspired by Thurston's conjecture.

\begin{question}
Let $M$ be a closed hyperbolic 3-manifold with positive first Betti number. Can every integral class $a \in H^2(M;\mathbb{R})$ (or $H^2(M ; \mathbb{Z})$) of dual Thurston norm exactly (at most) one (and satisfying the parity condition) be realized as the Euler class of some (faithful) representation $\rho \colon  \pi _1(M) \longrightarrow \Text{Homeo}^+(S^1)$?
\end{question}

Miyoshi gave a partial positive answer to the above question for closed orientable Seifert fibered manifolds, $H^2(M ; \mathbb{Z}) $, cohomology classes of dual Thurston norm at most one, and smooth representations \cite[Theorem 2]{miyoshi1997foliated}. Miyoshi also gave examples of cohomology classes that are representable as continuous Euler classes but not as smooth ones \cite[Theorem 1]{miyoshi1997foliated}; he showed, using a rigidity theorem of Ghys \cite{ghys1993rigidite}, that in fact the Euler class of any fibration of a closed hyperbolic 3-manifold is not smoothly representable. See also Miyoshi \cite{miyoshi2004representability}. 

Note that if $M$ is a closed aspherical 3-manifold, the Euler class of a representation $\rho \colon \pi_1(M) \rightarrow \text{Homeo}^+(S^1)$ does not necessarily satisfy the parity condition. For example, take $M = S_g \times S^1$ with $S_g$ a closed hyperbolic surface of genus $g$, and let $\rho$ be a representation of $ \pi_1(M) = \pi_1(S_g) \times \mathbb{Z}$ such that the image of $\mathbb{Z}$ is trivial, and the image of $\pi_1(S_g)$ has odd Euler class. See Calegari \cite[Remark 3.3]{calegari2006real}, and Culler and Dunfield \cite[Lemma 8.2]{culler2018orderability} for a parity condition for holonomy representations of cusped hyperbolic 3-manifolds. \\

\textbf{Virtual Euler class one conjecture}: Given a foliation $\mathcal{F}$ on $M$ with Euler class $a$ and a finite covering map $p \colon \tilde{M} \longrightarrow M$, the foliation $\mathcal{F}$ can be pulled back to a foliation $\tilde{\mathcal{F}}$ on $\tilde{M}$ with $e(\tilde{\mathcal{F}}) = p^*(a)$. By Proposition \ref{covering}, if $a \in H^2(M ; \mathbb{R})$ has dual norm equal to (respectively at most) one, then $p^*(a)$ has dual norm equal to (respectively at most) one.

\begin{question}
Let $a \in H^2(M; \mathbb{R})$ be an integral point with dual norm equal to (at most) one that satisfies the parity condition. Is there a finite covering map $p \colon \tilde{M} \longrightarrow M$ and a taut foliation $\tilde{\mathcal{F}}$ of $\tilde{M}$ such that $e(\tilde{\mathcal{F}})=p^*(a)$?
\end{question}

\bibliographystyle{plain}
\bibliography{foliationreference2}

\end{document}